\documentclass{amsart}
\let\tensor\otimes

\newtheorem{theorem}{Theorem}[section]
\newtheorem{conjecture}[theorem]{Conjecture}
\newtheorem{proposition}[theorem]{Proposition}
\newtheorem{lemma}[theorem]{Lemma}
\newtheorem{corollary}[theorem]{Corollary}
\theoremstyle{definition}
\newtheorem{definition}[theorem]{Definition}
\newtheorem{notation}[theorem]{Notation}
\newtheorem{remark}[theorem]{Remark}
\newtheorem{example}[theorem]{Example}
\newtheorem{problem}[theorem]{Problem}
\newtheorem{construction}[theorem]{Construction}
\theoremstyle{plain}
\numberwithin{equation}{section} 
\numberwithin{figure}{section} 

\def\({\left(}
\def\){\right)}
\def\a{\alpha}
\def\b{\beta}
\def\c{\gamma}
\def\C{{\mathbb C}}
\def\calC{{\mathcal C}}

\def\D{{\mathcal D}}
\def\dd{d}
\def\e{\varepsilon}
\def\F{{\mathbb F}}
\def\Gal{{\rm Gal}}
\def\hsin #1 {H_{#1}}
\def\INT{K_2(C;\Z)}
\def\iso{\cong}
\def\Ker{{\rm ker}}
\def\KCQbar{K_2^T(C_\Qbar)}
\def\O{{\mathcal O}}
\def\ord{\text{ord}}
\def\Res{{\rm Res}}
\def\s{\sigma}
\def\sheafK2{\underline{\mathcal K}_2}
\def\specelt{{\mathbb M}}
\def\tC{\widetilde C}
\def\tM{\widetilde M}
\def\torsion{{\rm torsion}}
\def\w{\wedge}

\catcode`\@=11

\def\newmathop#1{\expandafter\gdef\csname #1\endcsname{\mathop{\rm #1}\nolimits}}
\def\newmathopl#1{\expandafter\gdef\csname #1\endcsname{\mathop{\rm #1}\limits}}

\def\si#1{_{\scriptscriptstyle #1}}
\def\defsubscript#1{\expandafter\gdef\csname i#1\endcsname{_{\scriptscriptstyle #1}}}
\count0=26\loop\expandafter\expandafter\expandafter\defsubscript\@Alph{\the\count0}
\advance\count0 by -1\ifnum\count0>1\repeat

\def\defcalligraphic#1{\expandafter\gdef\csname c#1\endcsname
{{\ifmmode\mathcal#1\else$\mathcal#1$\fi}\spaceifletter}}
\count0=26\loop\expandafter\expandafter\expandafter\defcalligraphic\@Alph{\the\count0}
\advance\count0 by -1\ifnum\count0>1\repeat

\def\beq{\futurelet\comingchar\dobeq}     
  \def\dobeq{\ifx\comingchar*\begin{eqnarray*}\let\skipstar=\skiptoken
             \else\begin{eqnarray}\let\skipstar=\relax\fi\skipstar}
\def\eeq{\futurelet\comingchar\doeeq}     
  \def\doeeq{\ifx\comingchar*\end{eqnarray*}\let\skipstar=\skiptoken
             \else\end{eqnarray}\let\skipstar=\relax\fi\skipstar}
\def\skiptoken#1{}

\def\mdef#1#2{\def#1{\relax\ifmmode{#2}\else$#2$\fi\spaceifletter}}
\def\tdef#1#2{\def#1{\relax\ifmmode{\mathop{\rm #2}}\else#2\fi\spaceifletter}}

\def\spaceifletter{\futurelet\comingchar\dospaceifletter}
\def\dospaceifletter{\relax\ifmmode\else
  \ifcat A\noexpand\comingchar{} \fi
  \ifcat 0\noexpand\comingchar
  \ifx 0\noexpand\comingchar{} \fi
  \ifx 1\noexpand\comingchar{} \fi\ifx 2\noexpand\comingchar{} \fi
  \ifx 3\noexpand\comingchar{} \fi\ifx 4\noexpand\comingchar{} \fi
  \ifx 5\noexpand\comingchar{} \fi\ifx 6\noexpand\comingchar{} \fi
  \ifx 7\noexpand\comingchar{} \fi\ifx 8\noexpand\comingchar{} \fi
  \ifx 9\noexpand\comingchar{} \fi\fi
  \ifcat $\noexpand\comingchar{} \fi
  \ifcat \noexpand\relax\noexpand\comingchar{} \fi
\fi}

\def\isinlist#1#2#3{
\def\xisinlist##1##2#2##3##4\skiptoken{\ifx##3\skiptoken
  \def\isinlistfound{F}\else\def\isinlistfound{T}##1\fi
  }\xisinlist{#3}#1\donothing#2\skiptoken\skiptoken\donothing}

\def\inputdocument#1{\begingroup
  \def\documentstyle{\let\oldinput=\input\let\input=\skiptoken
  \let\inputifexists=\skiptoken
  \futurelet\comingchar\andocstyle}%
  \let\documentclass\documentstyle
  \def\andocstyle{\ifx\comingchar[\let\skipdoc=\skipdocstyle\else\let
  \skipdoc=\skiptoken\fi\skipdoc}%
  \def\document{\let\input=\oldinput\let\inputifexists=\oldinputifexists}%
  \def\enddocument{}%
  \input{#1}\endgroup}
\def\skipdocstyle[#1]#2{}

\catcode`\@=12

\mdef\P{\mathbb P}
\mdef\Q{\mathbb Q}
\mdef\R{\mathbb R}
\mdef\Z{\mathbb Z}
\mdef\C{\mathbb C}
\mdef\F{\mathbb F}

\def\iZ{\si\Z}

\let\tensor\otimes

\def\Qbar{{\overline\Q}}
\def\CZ{\cC}

\def\D{{\mathcal D}}

\mdef\KF{K_2(\Q(C))}                           
\mdef\KC{K_2^T(C)}                             
\mdef\KCQbar{K_2^T(C_\Qbar)}                   
\let\KCQB\KCQbar
\mdef\KCtors{K_2^T(C)/{\rm torsion}}           
\mdef\KCQ{K_2^T(C)\tensor\Q}                   

\newmathop{Jac}
\newmathop{div}
\newmathop{Div}
\newmathop{Spec}
\newmathop{rk}
\newmathop{Pic}
\newmathop{lcm}

\def\Definition#1\par{\begin{definition}#1\end{definition}\par}
\def\Conjecture#1\par{\begin{conjecture}#1\end{conjecture}\par}
\def\Construction#1\par{\begin{construction}#1\end{construction}\par}
\def\Example#1\par{\begin{example}#1\end{example}\par}
\def\Proposition#1\par{\begin{proposition}#1\end{proposition}\par}
\def\Theorem#1\par{\begin{theorem}#1\end{theorem}\par}
\def\Remark#1\par{\begin{remark}#1\end{remark}\par}
\def\Lemma#1\par{\begin{lemma}#1\end{lemma}\par}
\def\Notation#1\par{\begin{notation}#1\end{notation}\par}

\def\ben{\begin{enumerate}}
\def\een{\end{enumerate}}
\def\bit{\begin{itemize}}
\def\eit{\end{itemize}}

\def\ExampleNamed#1{\par\medskip\noindent\refstepcounter{theorem}{\bf Example \thetheorem: #1.}}

{\catcode`\@=11\global\let\c@equation=\c@theorem}

\input{xypic}
\usepackage{epsfig}
\begin{document}

\title[Beilinson's conjecture for $K_2$ of hyperelliptic curves]%
{Numerical verification of Beilinson's conjecture for $K_2$ of hyperelliptic curves.}
\author{Tim Dokchitser}
\address{Department of Mathematical Sciences\\University of Durham\\Science Laboratories\\South Road\\Durham DH1 3LE\\United Kingdom}
\author{Rob de Jeu}
\address{Department of Mathematical Sciences\\University of Durham\\Science Laboratories\\ South Road\\Durham DH1 3LE\\United Kingdom}
\author{Don Zagier}
\address{Max-Planck-Institut f\"ur Mathematik\\Vivatsgasse 7\\D-53111 Bonn\\Germany
{\rm and}
Coll\`ege de France\\ 11 place Marcelin Berthelot\\F-75005 Paris\\France}

\begin{abstract}
We construct families of hyperelliptic curves over $ \Q $
of arbitrary genus $ g $ with (at least) $ g $ integral elements
in $ K_2 $.  We also verify the Beilinson conjectures about $ K_2 $
numerically for several curves with $ g=2 $, $3$, $ 4 $ and $ 5 $.
The first few sections of the paper also provide an elementary introduction to the Beilinson conjectures
for $K_2$ of curves.
\end{abstract}

\maketitle

\section{Introduction} 
\label{intro}

Let $ k $ be a number field, with $ r_1 $ real embeddings and
$ 2r_2 $ complex embeddings into $ \C $, so that $ [ k : \Q ] = r_1 + 2 r_2 $.
It is a well known classical theorem that, if $ \O_k $ is the
ring of algebraic integers in $ k $, then
$ \O_k^* $ is a finitely generated abelian group of rank $ r = r_1 + r_2 - 1 $.
If $ u_1,\dots,u_r $ form a basis of $ \O_k^*/\torsion $, and $ \s_1,\dots,\s_{r+1} $
are the complex embeddings of $ k $ up to complex conjugation,
then the {\it regulator\/} of $ \O_k^* $ is defined by
\[
R = \dfrac{2^{r_2}}{[k:\Q]}  \left | 
\det
\left(
\begin{matrix}
1 & \log|\s_1(u_1)| & \dots & \log|\s_1(u_r)| \cr
\vdots & \vdots & & \vdots \cr
1 & \log|\s_{r+1}(u_1)| & \dots & \log|\s_{r+1}(u_r)|
\end{matrix}
\right)
\right |  
\; ,
\]
and one has that
$ \Res_{s=1}\zeta_k(s) =  \dfrac{2^{ r_1 } (2\pi)^{ r_2 } R \, |\text{Cl}(\O_k)| }{ w \sqrt{  \Delta_k } } $,
where $ \Delta_k $ is the absolute value of the discriminant of $ k $
and $ w $ the number of roots of unity in $ k $.
Because  $ K_0(\O_k) \iso \text{Cl}(\O_k) \oplus \Z $ and $ K_1(\O_k) \iso \O_k^* $, 
so $ |\text{Cl}(\O_k)| = |K_0(\O_k)_{\text{tor}}| $
and $ w = | K_1(\O_k)_{\text{tor}} | $, this can be interpreted
as a statement about the
$ K $-theory of $ \O_k $, and it is from this point of view that
it can be generalized to $ \zeta_k(n) $ for $ n \geq 2 $.
Namely, in \cite{Quillen73}, Quillen proved that $ K_n(\O_k) $ is a finitely
generated abelian group for all~$ n $.
Borel in \cite{Borel74} computed its rank, showing that this
rank is zero for even $ n \geq 2 $ and is equal to $ r_\pm $ for
odd $ n = 2m -1 > 1 $,
where $ (-1)^m = \pm 1 $ and $ r_- = r_1 + r_2 $, $ r_+ = r_2 $.
Moreover, for those odd $ n $
he showed (see \cite{Borel77}) that
a suitably defined regulator of $ K_{2m-1}(\O_k) $ 
is a non-zero rational multiple of
$ \zeta_k(m)/\pi^{m r_\mp} \sqrt{\Delta_k} $.

Inspired by this, Bloch in \cite{Bloch} considered $ K_2 $ of elliptic
curves $ E $ defined over $ \Q $ with complex multiplication, and showed
that there is a relation between a regulator associated to
$ K_2(E) $ and the value of $ L(E, 2) $. 
Beilinson then proposed a very general conjecture about similar relations between certain regulators
of $ K $-groups of projective varieties over number fields and values of their
$ L $-functions at integers (see \cite[\S5]{Schn}).

Those conjectures were tested numerically for $ K_2 $ of elliptic curves over
$ \Q $ by Bloch and Grayson in  \cite{BlGr}, which led to a modification
of Beilinson's original conjecture, as explained in Section~3.  Some further numerical work has been
done in this direction, e.g., 
Young carried
out similar calculations for elliptic curves over certain real quadratic
number fields (as well as over $ \Q $) in his thesis \cite{young},
and 
Kimura worked out the case of the genus two curve 
$ y^2 - y = x^5 $ in~\cite{Kim}.

The goal of this paper is to verify Beilinson's conjecture numerically
for $ K_2 $ of a number of hyperelliptic curves over $ \Q $ of genus greater 
than 1.
The results of the computations support the predictions of Beilinson's conjecture in the cases we studied.
Specifically, the conjecture says that the rank of a certain torsion-free abelian group $ \INT $
associated to the curve $ C $ over $ \Q $ is equal to the genus $ g $ of $ C $, and that the associated
regulator (the determinant of a certain $ g \times g $-matrix) is rationally proportional to the appropriately
normalized value of $ L(C,2) $ (see~Section~3).
If $ g>1 $ it is quite difficult to write down enough elements in $ \INT $, but we give
several infinite families of hyperelliptic curves of genus 2 and 3 and one further infinite family for every genus $ g \geq 2 $,
as well as sporadic examples for $ g=4 $ or 5,
for which we can construct at least $ g $ elements of $ \INT $.  For some 200 of those curves we check
by computer that the regulator of $ g $ of those elements is non-zero and is related to the $ L $-function of the
curve in the expected way.
Since the verification that a real number is non-zero can be done numerically, these calculations prove rigorously that
the $ g $ elements in question are linearly independent and hence that $ \rk \INT \geq g $ for these
curves.
The relationship between the regulator and $ L(C,2) $, on the other hand, can only be established numerically to high precision.

Furthermore, the second author has shown, in a sequel to this paper (\cite{dJ}), that in the family for
arbitrary $ g \geq 2 $ there are, for each $ g $, infinitely many non-isomorphic curves for which the regulator of the
$ g $ elements is non-zero, thus showing that $ \rk \INT \geq g $ for those curves.  For a more precise formulation
of those (and other) results we refer the reader to Remark~\ref{sequel} or \cite{dJ}.

Our results also provide evidence
for the reverse inequality $ \rk \INT \leq g $ predicted by Beilinson's conjecture.  Namely, in a number of sporadic
examples, and one universal situation (see Remark~\ref{seven} for the precise list), we actually construct
more than $ g $ elements of $ \INT $.  In each case the computer calculations of the regulator
suggested a linear dependence over $ \Z $ for any $ g+1 $ of our elements.  In most cases, including
the universal case, these relations could then be proved;  in the few remaining cases we have only a
high-precision numerical verification at the regulator level.  Since there is no intrinsic reason why
the elements we construct should be linearly dependent when there are more than $ g $ of them, these results can be
seen as evidence for the prediction that the rank of $ \INT $ is at most $ g $.

The structure of the paper is as follows.
In Sections~2 and~3 we review the statement of Beilinson's conjecture for the
case of $ K_2 $ of curves defined over $ \Q $.
In the next four sections we show how to construct curves
with interesting elements in $ K_2 $.
In Section~8 we deal with a technical condition in Beilinson's
conjecture, the {\it integrality condition\/}, and in Section~9
we discuss how to compute the regulator numerically.  The final
section is devoted to examples.

The first two authors would like to thank the Coll\`ege de France for
a visit to work on this paper, as well as the EC network Arithmetic
Algebraic Geometry for travel support.  The second
author would also like to thank the Newton Institute for a productive stay during
which this paper was worked on.

\section{Curves and their $ L $-functions  (review)}  
\label{lsection}

Let $ C $ be a non-singular, projective, geometrically irreducible curve over $ \Q $ of genus $ g $.  
Its $ L $-function is the Dirichlet series defined for $ \text{Re}(s) > \frac{3}{2} $ by the Euler product
$$
L(C,s) = \prod_{\text{$ p $ prime}} L_p(C,s) 
\; ,
$$
where the Euler factor $L_p(C,s)$, for those primes $p$ for which the equations
defining $ C $ over $ \Q $ can be reduced modulo $ p $ and still define a
smooth curve of genus $ g $ over $ \F_p $, is defined by the formula
$$
L_p(C,s) = \exp\bigg( \sum_{n=1}^\infty \big( p^n+1-\#C(\F_{p^n}) \big) \frac{p^{-ns}}{n} \bigg)
$$
and is known to be the reciprocal of a polynomial in $ p^{-s} $ of degree
$ 2g $ with constant term 1.  For the remaining (finitely many) primes $ p $,
$ L_p(C,s)^{-1} $ is also a polynomial in $ p^{-s} $ with constant
term 1, but now
of degree at most $ 2g $.  
Finally, one can associate to the curve  $  C/\Q $ a positive
integer $ N $ called its {\it conductor} which plays a role in the
expected functional equation of its $L$-function (see below).  We omit the precise
definitions of the remaining Euler factors and of $ N $ since 
they are a little complicated and since these quantities
can be (and in some cases actually were) computed experimentally by assuming the functional
equation rather than from their definitions.  This functional equation is:

\begin{conjecture}{\bf (Hasse-Weil)}
\label{lconjecture}
The function
$$
L^*(C,s) = \frac{N^{s/2}}{(2 \pi)^{gs}} \, \Gamma(s)^g  \, L(C,s)
$$
extends to an entire function of $ s $ and
satisfies 
$ L^*(C,s) = w L^*(C,2-s)  $, where $ w = + 1 $ or $ - 1 $.
\end{conjecture}

Note that, if this conjecture holds, then 
$ L^{(0)}(C,0) = \dots = L^{(g-1)}(C,0) = 0 $ and
\begin{equation*}
\frac{L^{(g)}(C,0)}{g!} =
\lim_{s \to 0} \frac{L(C,s)}{s^g} =
L^*(C,0) =
w L^*(C,2) =
\frac{wN}{(2 \pi)^{2g}}  L(C,2) \neq 0
\; ,
\end{equation*}
because $ \Gamma(s) $ has a pole of order one at $ s=0 $ with
residue 1.

\begin{remark}
In this paper we need the value of $ L(C,s) $
only at $ s=2 $, where the defining Euler product is absolutely
convergent.  However, this convergence is very slow and in practice
we will calculate the value of $ L(C,2) $ numerically later
on by assuming that Conjecture~\ref{lconjecture} holds and 
using the algorithms described in \cite{Dok1} and~\cite{Dok2}.
These algorithms also compute $ N $, $ w $, and the Euler factors
$ L_p(C,s) $ for \lq\lq bad\rq\rq\ primes
under the assumption that the functional equation of $ L^*(C,s) $ holds for some choices
of these quantities from a finite list of possibilities,
and at the same time check the validity of the resulting functional equation for $ L^*(C,s) $ to high
precision.
\end{remark}

\section{$K$-theory, regulators, and the Beilinson conjectures} 

In this section we will give definitions of the various $ K $-groups
occurring that are more elementary and explicit than the usual
ones, but that are equivalent in our situation.  We will indicate
the relations to the standard definitions as we go along.

Let $ F $ be a field.   Then $ K_2(F) $ can be 
defined\footnote{The actual definition of  $K_2$ (of arbitrary rings) is more complicated
(see~\cite[\S5]{Milnor}), and its equivalence for fields with the definition in
terms of the \lq\lq symbols\rq\rq\ $ \{a,b\} $ is a famous theorem of Matsumoto
(see \cite[Theorem~11.1]{Milnor}).}
as
$$
F^* \tensor_\Z F^* / \langle a \tensor (1-a), a \in F, a \neq 0,1 \rangle 
\; ,
$$
where $ \langle\cdots \rangle $ denotes the subgroup
generated by the elements indicated.
The class of $ a \tensor b $ is denoted $ \{ a, b \} $, so
that $ K_2(F) $ is an abelian group (written additively), with
generators $ \{ a , b \} $ for $ a $ and $ b $ in $ F^* $, and relations
\begin{alignat*}{1}
& \{ a_1a_2 ,  b \} = \{ a_1 ,  b \} + \{ a_2 ,  b \}
\\
& \{ a , b_1b_2 \} = \{ a , b_1 \} + \{ a , b_2 \}
\\
& \{ a , 1-a \} = 0 \text{ if } a \text{ is in } F , \ a \neq 0, 1 .
\end{alignat*}
It is a nice exercise to show that these relations imply that
$ \{a,b\} = - \{b,a\} $ and $  \{ c , -c \} = 0 $ for $ a $,
$ b $ and $ c $ in $ F^* $.

Now consider a (non-singular, projective, geometrically irreducible) curve $ C $ defined over $ \Q $.
Associated to $ C $ are $ K $-groups $ K_n(C) $ whose definition
is a little complicated, but for this paper we need only a certain
quotient group\footnote{usually denoted $ H^0(C,\sheafK2) $}
$ K_2^T(C) $ of $ K_2(C) $ (\lq\lq T\rq\rq\ for \lq\lq tame\rq\rq), which can be described in a simpler way.
Set
\begin{equation*}
K_2^T(C) = \Ker\bigg( K_2(F) \overset T \longrightarrow \bigoplus_{ x \in C(\Qbar) } \Qbar^* \bigg)
\; ,
\end{equation*}
where $ F = \Q(C) $ is the function field of $ C $
and where the $ x $-component of the map $ T $ is the {\it tame symbol\/}
at $ x $, defined by
\begin{equation}
\label{tame}
T_x : \{ a , b \} \mapsto
 (-1)^{ \ord_x(a)\ord_x(b)} \, \dfrac{a^{\ord_x(b)}}{b^{\ord_x(a)}} (x)
\; .
\end{equation}
Note that this definition makes sense
since $ \dfrac{a^{\ord_x(b)}}{b^{\ord_x(a)}} $
has order zero at $ x $ and hence is defined and non-zero at $ x $.
It is clear that $ T_x $ is a map on $ F^* \tensor_\Z F^* $, and
checking that it is trivial on the symbols $ \{ a, 1-a \} $ for $ a \neq 0,1 $
is a good exercise (which also explains why we want to have the
power of $ -1 $ in the formula), so $ T_x $ defines a map on
$ K_2(F) $.

We note that if $ \a $ is an element of $ K_2(\Q(C)) $,
then 
\begin{equation}\label{productformula}
\prod_{x \in C(\Qbar)} T_x(\a) = 1
\; ,
\end{equation}
a result known as the {\it product formula} (see~\cite[Theorem~8.2]{Bass}).

Beilinson, generalizing work by Bloch \cite{Bloch},
defined regulators of the $ K $-groups of $ C $
(see \cite{Schn}).
We will describe these in elementary terms for $ K_2^T(C) $.

We start with a map from $ F^* \times F^* $ to the group of almost
everywhere defined 1-forms on the Riemann surface $ X = C(\C) $
by putting
\begin{equation}
\label{eta}
\eta( a , b) =  \log|a| \, \dd \arg b - \log|b| \, \dd \arg a 
\; ,
\end{equation}
where $ \arg a $ is the argument of $ a $.  Note that this  is
well defined (the argument $ \arg $ is defined up to multiples of
$ 2 \pi $, but these map to zero under $ \dd $) and is
a smooth (indeed, real-analytic) 1-form on the complement of the
set of zeros and poles of $ a $ and $ b $.
It is clear that $ \eta( a_1 a_2 , b ) = \eta(a_1,b) + \eta(a_2,b) $
and $ \eta(a,b_1 b_2) = \eta(a, b_1) + \eta(a, b_2) $, so that
$ \eta $ induces a map (still denoted $ \eta $) on $ F^* \tensor_\Z F^* $.
Also, $ \eta(a,b) $ is closed, since $ \dd \eta(a , b) = \text{Im} ( \dd \log a \w \dd \log b ) = 0 $.

For any smooth closed 1-form $ \omega $ defined on the complement
of a finite set $ S \subset X $, and any smooth oriented loop
$ \c $ in $ X \setminus S $, we have a pairing
\[
 ( \gamma,\omega ) = \frac{1}{2\pi} \int_\gamma \omega
\]
which depends only on the homology class of $ \c $ in $ X \setminus S $.
As $ \c $ moves across a point $ x $ in $ S $, the value of 
$ ( \c, \omega ) $ jumps by $ ( C_x , \omega )  $,
where $ C_x $ denotes a small circle around $ x $.  A simple
calculation shows that $ ( C_x , \eta(a,b) ) = \log | T_x(\{a,b\})| $.
It follows that if $ \a = \sum_i a_i \tensor b_i $ is an element
of $ F^* \tensor_\Z F^* $ such that $ T_x(\a) = 1 $ for all $ x $
in $ X $, then $ ( \, \cdot \, , \eta(\a) ) $ is a well
defined map from $ \hsin 1 (X ; \Z ) $ to $ \R $.

Next, one has to check that this map vanishes if $ \a = a \tensor (1-a) $ and hence gives us a pairing 
\[
\langle \, \cdot \, , \, \cdot \, \rangle : \hsin 1 (X ; \Z) \times K_2^T(C) \to \R
\]
given by $ \langle \gamma, \a \rangle = ( \gamma, \eta(\a) ) $.
This follows from the fact that $ \eta(a,1-a)=\dd\,D(a) $, where
$ D(z) $ is the Bloch-Wigner dilogarithm function.

Finally, we observe that, if $ c $ denotes complex conjugation
on $ X $, then $ c^*(\eta(\a)) = - \eta(\a) $ for any $ \a $ in $ K_2^T(C) $,
since $ c^*(\log|a|) = \log|a| $ and $ c^* ( \dd \arg b ) = - \dd \arg b $:
$ a $ and $ b $ are in $ \Q(C)^* $, so $ (c^* b)(z) = b(\overline z )  = \overline {b(z)}$.
This means that if $ \gamma $ is in $ \hsin 1 (X ; \Z)^+ $, the
$ c $-invariant part of $ \hsin 1 (X ; \Z) $,
then $ \langle \gamma,\a \rangle = 0 $ for any $ \a $ in $ K_2^T(C) $:
\[
\langle \gamma,\a \rangle =
\frac{1}{2\pi}\int_\gamma \eta(\a) =
\frac{1}{2\pi}\int_{c \circ \gamma } \eta(\a) =
\frac{1}{2\pi}\int_\gamma c^*(\eta(\a)) =
- \frac{1}{2\pi}\int_\gamma \eta(\a) =
-\langle \gamma,\a \rangle
\; .
\]
We therefore have to compute $ \langle \gamma,\a \rangle $ only
for $ \gamma $ in
$ \hsin 1 (X; \Z) / \hsin 1 (X; \Z)^+ $.  In practice,
we can just as well compute it for all $ \gamma $ in $ \hsin 1 (X ; \Z)^- $,
the anti-invariants in $ \hsin 1 (X ; \Z ) $ under the action
of $ c $, giving us finally the {\it regulator pairing\/}
\begin{equation}
\label{regulatorpairing}
\langle \, \cdot \,  , \, \cdot \, \rangle : \hsin 1 (X ; \Z)^- \times K_2^T(C)/\torsion \to \R
\; .
\end{equation}

It is easy to see that $ \hsin 1 (X ; \Z)^- $
has rank $ g $.
Beilinson originally conjectured\footnote{Strictly
speaking, Beilinson conjectured
$ K_2^T(C) \tensor_\Z \Q $ to have dimension $ g $ over $ \Q $.
However, $ K_2^T(C)/K_2(\Q) $ is expected to be finitely generated
by virtue of a conjecture of Bass, and since $ K_2(\Q) $ is a
torsion group, this would imply the equivalence of Beilinson's
formulation and our formulation.}
that the rank of $ K_2^T(C)/\torsion $
is also equal to $ g $, that the pairing in \eqref{regulatorpairing} was non-degenerate,
and that there was a relation between $ L(C,2) $ and the absolute
value of the determinant of the matrix of this pairing with
respect to bases of $ \hsin 1 (X, \Z)^- $ and $ K_2^T(C)/\torsion $.
Unfortunately, this conjecture was 
wrong, as $ K_2^T(C)/\torsion $ can have rank
bigger than $ g $ already for $ g=1 $, as was discovered by
Bloch and Grayson in \cite{BlGr}.
They found that one should consider a certain subgroup of $ K_2^T(C)/\torsion $
defined by an additional condition, which we now proceed to
describe.  

The extra condition comes from the fact that one should not consider
the curve over $ \Q $, but instead a model of it over $ \Z $.
This is analogous to the situation in Section~\ref{intro}, where
one has to consider $ K_1(\O_k) \iso \O_k^* $ instead of $ K_1(k) \iso k^* $
in order to get the correct regulator.

So let $ \calC $ be a regular proper model of $ C $ over $ \Z $,
i.e., a regular, proper, irreducible two-dimensional scheme over $ \Z $ 
such that the generic fiber $ \calC_\Q $ is isomorphic
to $ C $ (see, e.g.,  \cite[Chapter~10]{Liu}).
For each prime $ p $, let $ \calC_p $ be the fiber
of $ \calC $ over $ \F_p $.  For each irreducible component $ \D $
of the curve $ \calC_p $, let $ \F_p(\D) $ denote its field of rational
functions over $ \F_p $.  Then we define
\begin{equation}
\label{KTcalC}
K_2^T(\calC) = \ker \bigg( K_2^T(C) \to \bigoplus_{p\!,\, \D \subseteq \calC_p}  \F_p(\D)^* \bigg)
\; ,
\end{equation}
where the map to $ \F_p(\D)^* $ is given as follows.  
The order of vanishing
along $ \D $ gives rise to a discrete valuation on $ F $, $ v_\D $.
The component of the map in \eqref{KTcalC} corresponding to $ \D $ is given by
the tame symbol corresponding to $ \D $,
\begin{equation}
\label{tameD}
T_\D : \{ a , b \} \mapsto
 (-1)^{ v_\D(a)v_\D(b)} \dfrac{a^{v_\D(b)}}{b^{v_\D(a)}} (\D)
\; ,
\end{equation}
in complete analogy with~\eqref{tame}.
Finally, we set 
\begin{equation}
\label{INTdef}
\INT = K_2^T(\calC)/\torsion
\> ,
\end{equation}
a subgroup of $ K_2^T(C)/\torsion $.  This group, sometimes denoted by
$ K_2(C)_\Z/\torsion $ (cf.~\cite[\S~3.3]{Soule}), is independent of the choice of the
regular proper model $ \calC $ of $ C $ (see~\cite[page~13]{Schn}),
justifying the notation.

\begin{remark}
We could have defined $ K_2^T(\calC) $ in a single step as
\begin{equation}
\label{realtame}
K_2^T(\calC) = \ker \bigg( K_2(F) \to \bigoplus_{\D} \F(\D)^* \bigg)
\; ,
\end{equation}
where $ \D $ runs through all irreducible curves on $ \calC $
and $ \F(\D) $ stands for the residue field at $ \D $.
Any such $ \D $ is either \lq\lq vertical\rq\rq, in which case
it is a component of some $ \calC_p $ and $ T_\D $ is the map
in~\eqref{tameD}, or else \lq\lq horizontal\rq\rq, in which case
it corresponds to the $ \Gal(\Qbar/\Q) $-orbit of some $ x $
in $ C(\Qbar) $ and $ T_\D $ being trivial is equivalent to
$ T_y $ being trivial for all $ y $ in that $ \Gal(\Qbar/\Q) $-orbit.
\end{remark}

We can now restrict the pairing \eqref{regulatorpairing} to
\begin{equation}
\label{INTpairing}
\langle \, \cdot \,  , \, \cdot \, \rangle : \hsin 1 (X ; \Z)^- \times \INT \to \R
\; ,
\end{equation}
and formulate our description of Beilinson's conjecture, as modified in
accordance with \cite{BlGr}, as follows:

\begin{conjecture}
\label{beilinsonconjecture}
Let $ C $ be a non-singular, projective, geometrically irreducible
curve of genus $ g $ defined over $ \Q $, and let $ X = C(\C) $.  Then:
\begin{enumerate}
\item
The group $ \INT $ is a free abelian group of rank $ g $
and the pairing~\eqref{INTpairing}
is non-degenerate;
\item
Let $ R $ denote the absolute value of the determinant of this
pairing with respect to $ \Z $-bases of $ \hsin 1 (X; \Z)^- $
and $ \INT $, and let $ L^*(C,0) $ be defined as in
Section~\ref{lsection}.  Then $ L^*(C,0) = Q R $
for some non-zero rational number~$ Q $.
\end{enumerate}
\end{conjecture}

\begin{remark}
The definition of $ L^*(C,0) $ requires the analytic continuation
of $ L(C,s) $, but since the analytic continuation and the expected
functional equation of $ L(C,s) $ would imply that $ L^*(C,0) $ is
rationally proportional to $ \pi^{-2g} L(C,2) $, Beilinson's
conjecture could be formulated without any assumptions about the
analytic continuation of $ L(C,s) $.
\end{remark}

\begin{remark}
\label{untractable}
In practice, the conjecture is rather intractable,
as it seems impossible to compute $ \INT $ even
after tensoring this group with $ \Q $.  Indeed, we neither
can guarantee finding enough elements to generate $ \INT $,
nor can we necessarily determine the rank of a subgroup generated by finitely
many elements, as we do not know any practical method for determining
if a given combination of elements in $ F^* \tensor _\Z F^* $
can be written as a sum of Steinberg symbols $ a \tensor (1-a) $.
But we can try to find $ g $ elements in $ \INT $ and
compute $ R $ as in the conjecture using those elements rather
than a basis of $ \INT $.  If $ R $ is non-zero numerically
we can check the relation with $ L^*(C,0) $ as in Conjecture~\ref{beilinsonconjecture}.
Also, if we have more than $ g $ elements $ \a_j $ in $ \INT $,
the conjecture implies that the maps
$ \langle \, \cdot \, , \a_j \rangle : \hsin 1 (X ; \Z)^- \to \R $
should be linearly dependent over $ \Z $, and this too can be checked numerically.
Both types of verification will be carried out in Section~\ref{examples}.
\end{remark}

\begin{remark}
We have restricted ourselves to the statement of Beilinson's
conjecture for a curve over $ \Q $, but the conjecture can be
formulated equally well for any number field $ k $.  Suppose that $ C/k $
is a (non-singular, projective, geometrically irreducible) curve of genus $ g $.  
Let $ \O_k $ be the ring of integers of $ k $, and let $ \calC $
be a model of $ C/k $ over $ \O_k $.  Then one defines $ K_2^T(\calC) $
as in \eqref{realtame}, with the sum over  irreducible curves
$ \D $ in $ \calC $.  Again, for the \lq\lq horizontal\rq\rq\
curves $ \D $, the $ T_\D $ correspond to the $ T_x $ for $ x $
in $ C(\Qbar) $, up to conjugation under $ \Gal(\Qbar/k) $.
Once again, $ K_2^T(\calC)/\torsion $ is independent of the choice of
the model $ \calC $, and is denoted by $ \INT $.
One expects $ \INT \iso \Z^{g \, [k:\Q]} $.
The main difference is that the Riemann surface involved will
no longer be connected.  Instead, let $ X $ be the Riemann
surface obtained from all points in $ C $ over $ \C $, using
all embeddings of $ k $ into $ \C $.  This is a disjoint
union of $ [k:\Q] $ connected Riemann surfaces, each of genus $ g $.
Complex conjugation acts on this either by swapping conjugate
pairs of complex embeddings of $ k $, or by acting on the Riemann
surface corresponding to a real embedding of $ k $.  Then 
$ \hsin 1 (X; \Z)^- \iso \Z^{g \, [k : \Q]} $, and there is
a pairing $  \hsin 1 (X; \Z)^- \times \INT \to \R $ given by
$ \langle \c , \a \rangle = \frac{1}{2\pi} \int_\c \eta(\a) $.
Here $ \eta(\a) $ is the 1-form on $ X $
which for $ \a = \{a , b\} $ is given 
on the connected Riemann surface corresponding
to $ \s : k \to \C $ by
$ \log|a_\s| \dd \arg b_\s - \log|b_\s| \dd \arg a_\s $.
(The subscripts indicate that we consider the functions 
on the Riemann surface obtained by applying $ \s $ to the coefficients
involved in $ a $ and $ b $.)
One again expects $ \langle \, \cdot \, , \, \cdot \, \rangle $
to be non-degenerate, defines $ R $ to be the absolute value of the determinant of the matrix
of $ \langle \, \cdot \, , \, \cdot \, \rangle $ with respect to $ \Z $-bases of
$ \hsin 1 (X ; \Z)^-  $ and $ \INT $, and conjectures that $ L^*(C,0) = Q R $
for some non-zero rational number~$ Q $.
\end{remark}


\section {Constructing elements of $K_2$ from torsion divisors} 
\label{sTorssymb}

The first problem in testing Conjecture~\ref{beilinsonconjecture} is that it is not
at all clear how to construct elements of $\KC$ on a given curve $C$
over $ \Q $,
that is, how to produce rational functions $f_i$, $g_i$ on $C$ such that
$\sum_i\{f_i,g_i\}$ has trivial tame symbol at every point of $C$.

To understand this condition, consider one symbol $\{f,g\}$ in $\KF$.
If $\div(f)$ and $\div(g)$ have disjoint support, then this symbol
lies in $\KC$ if and only if
$f(P)^{\ord\iP(g)}\!=\!1$ for every zero or pole $P$ of $g$ and
$g(Q)^{\ord\si Q(f)}\!=\!1$ for every zero or pole $Q$ of $f$. Essentially
this says that $f$ equals one (or a root of unity, if $|\ord\iP(g)|>1$)
whenever $g$ has a zero or a pole, and similarly with $ f $ and $ g $ interchanged.

To try to satisfy these conditions, it is natural to look at functions
which have only very few zeros and poles. The simplest case is given by
functions $f$ and $g$ which have only one (multiple) zero and one
(multiple) pole. If one of these points is common for $f$ and $g$,
then it turns out that simply renormalizing the functions is
sufficient to satisfy the tame symbol conditions. For this we use the
product formula~\eqref{productformula}.
All of our examples are then based on the following construction.

\Construction\label{constors}
Let $C/\Q$ be a curve. Assume $P_1, P_2, P_3 \!\in\! C(\Q)$ are distinct
points whose pairwise differences are torsion divisors. Thus there are rational
functions $f_i$ with
\beq*
  \div (f_i) = m_i (P_{i+1})-m_i(P_{i-1}), \qquad  i \in \Z/3\Z \>,
\eeq*%
where $m_i$ is the order of $(P_{i+1})\!-\!(P_{i-1})$ in the divisor group
$\Pic^0(C)$.  We then define three elements of \KF by
$$
  S_i = \Bigl\{ \frac{f_{i+1}}{f_{i+1}(P_{i+1})}, \frac{f_{i-1}}{f_{i-1}(P_{i-1})} \Bigr\}, \qquad i \in \Z/3\Z \>.
$$
The functions $f_i$ are unique up to constants, so the symbols $S_i$ are uniquely
defined by the points $P_i$.  Moreover, they satisfy the tame symbol condition everywhere:

\Lemma The $S_i$ are elements of $\KC$.
\begin{proof}
The components of $S_i$ are normalized
to make the tame symbol trivial at $P_{i-1}$ and $P_{i+1}$.
By the product formula (\ref{productformula})
it is also trivial at $P_i$, this being the only other
point in the support of the divisors of $f_{i-1}$ and $f_{i+1}$.
\end{proof}

\noindent
Next we show that the three elements $S_i$
generate a subgroup of rank at most 1 of~$\KCtors$.

\Proposition\label{consprop} We keep the notation of Construction \ref{constors}.
\begin{enumerate}
\item There is a unique element $\{P_1,P_2,P_3\}$ of $\KCtors$ such that
in this group
$$
   S_i \>\>=\>\> {\frac{\lcm(m_1,m_2,m_3)}{m_i}}\> \{P_1,P_2,P_3\}, \qquad i=1,2,3\>.
$$
\item The element $\{P_1,P_2,P_3\}$ is unchanged under even permutations
and changes sign under odd permutations of the points.
\end{enumerate}

\begin{proof}
(1) 
Uniqueness is obvious, so we only need to show existence.
Replacing the functions $f_i$ by $f_i/f_i(P_i)$ if necessary, we can
assume that $f_i(P_i)\!=\!1$, so that $ S_i = \{f_{i+1},f_i \} $.
Now let
\begin{equation}
\label{L}
  L = \lcm(m_1,m_2) = \lcm(m_1,m_3) = \lcm(m_2,m_3)\>;
\end{equation}
the three least common multiples are equal since every integer which kills both
$(P_i)\!-\!(P_{i+1})$ and $(P_{i+1})\!-\!(P_{i+2})$ in $\Pic^0(C)$ also kills
$(P_{i})\!-\!(P_{i+2})$,
which is their sum. It follows that $ L=\lcm(m_1,m_2,m_3) $ and that the quotients
$$
  r_i = L/m_i, \qquad i=1,2,3
$$
are pairwise relatively prime. Then
$$
  \div (f_i^{r_i}) = L(P_{i+1})\!-\!L(P_{i-1}),
$$
so the function $f_1^{r_1}f_2^{r_2}f_3^{r_3}$ has trivial divisor,
and is therefore a constant:
\beq\label{prodfc}
  f_1^{r_1}f_2^{r_2}f_3^{r_3} = c \>.
\eeq
Next, we have, in \KC
$$
  r_1 S_2 = \{f_3, f_1^{r_1}\}, \qquad
  r_2 S_1 = \{f_2^{r_2}, f_3\} = -\{f_3, f_2^{r_2} \}\>,
$$
so
$$
  r_1 S_2 - r_2 S_1 = \{ f_3, f_1^{r_1}f_2^{r_2} \}
    = \{ f_3, c f_3^{-r_3} \} = \{f_3,(-1)^{r_3}c\} \>,
$$
because $\{f_3, -f_3\}$ is trivial.
Since the left-hand side is in \KC, so is $\{f_3,(-1)^{r_3}c\}$.
This implies that $c^{2m_3}\!=\!1$, so $\{f_3,(-1)^{r_3}c\}$ is torsion.
Therefore, modulo torsion,
for $i,j=1,2,3$,
$$
  r_i S_j \!=\! r_j S_i \>.
$$
Now choose integers $\alpha_1, \alpha_2, \alpha_3$ such that $\sum \alpha_i r_i \!=\! 1$
and set
$$
  T = \alpha_1 S_1 + \alpha_2 S_2 + \alpha_3 S_3 \>.
$$
Then
$$
  r_1 T = \alpha_1 r_1 S_1 + \alpha_2 r_1 S_2 + \alpha_3 r_1 S_3 =
          \alpha_1 r_1 S_1 + \alpha_2 r_2 S_1 + \alpha_3 r_3 S_1 = S_1,
$$
and, similarly, $r_2 T \!=\! S_2$ and $r_3 T \!=\! S_3$, all
modulo torsion. In particular, the subgroup $\langle S_1, S_2, S_3 \rangle$
of $\KCtors$ is generated by the class of $T$ alone.


(2)
The first statement follows from the construction in part (1). For the second,
let $f_i$ be as above, chosen so that $f_i(P_i)\!=\!1$.

Take $\tilde P_1\!=\!P_2, \tilde P_2\!=\!P_1, \tilde P_3\!=\!P_3$
and the functions
$\tilde f_1\!=\!f_2^{-1}, \tilde f_2\!=\!f_1^{-1}, \tilde f_3\!=\!f_3^{-1}$,
so that $ \tilde\alpha_1 = \alpha_2 $, $ \tilde\alpha_2 = \alpha_1 $
and $ \tilde\alpha_3 = \alpha_3 $.
The corresponding symbols are
$$
  \tilde S_1 = \{\tilde f_2, \tilde f_3\} = \{f_1^{-1}, f_3^{-1}\} =
    - \{f_1, f_3^{-1}\} = \{f_1, f_3\} = -\{f_3, f_1\} = - S_2,
$$
and, similarly, $\tilde S_2\!=\!\{f_3^{-1},f_2^{-1}\}\!=\!-S_1$
and $\tilde S_3\!=\!\{f_2^{-1},f_1^{-1}\}\!=\!-S_3$.  Then modulo torsion,
$$
  \tilde T = \tilde\alpha_1 \tilde S_1 + \tilde\alpha_2 \tilde S_2 + \tilde\alpha_3 \tilde S_3 =
  - \alpha_2 S_2 - \alpha_1 S_1 - \alpha_3 S_3 = - T \>.
$$
\end{proof}

\Proposition\label{tetrel}
Let $C/\Q$ be a curve and $P_1, P_2, P_3, P_4$ be four distinct
points in $C(\Q)$ such that all $(P_i)\!-\!(P_j)$ are torsion divisors.
Then the four elements
$$
  \{P_1, \ldots, \widehat{P_i}, \ldots, P_4\} \qquad (1\le i\le 4)
$$
are linearly dependent. More precisely, if
$m_{ij}$ is the order of $(P_i)\!-\!(P_j)$, then
\begin{equation}
\label{Crelation}
  \sum_{i=1}^4 (-1)^i\> c_i \> \{P_1, \ldots, \widehat{P_i}, \ldots, P_4\} = 0
\end{equation}
holds in $ \KC/\torsion $ with
\begin{equation}
\label{Cis}
   c_i = \gcd(m_{ij}, m_{ik}, m_{il})
  \qquad
  (\{i,j,k,l\}\!=\!\{1,2,3,4\})
  \>.
\end{equation}

\begin{proof}
Choose $f_{ij}$ with $\div (f_{ij})\!=\!m_{ij}(P_i)\!-\!m_{ij}(P_j)$.
The discussion surrounding (\ref{prodfc}) shows that
$$
  \Bigl(\frac{f_{ij}}{f_{ij}(P_k)}\Bigr)^{\frac M{m_{ij}}}
  \Bigl(\frac{f_{jk}}{f_{jk}(P_i)}\Bigr)^{\frac M{m_{jk}}}
  \Bigl(\frac{f_{ki}}{f_{ki}(P_j)}\Bigr)^{\frac M{m_{ki}}} = 1
$$
if $\{i,j,k\}\!\subset\!\{1,2,3,4\}$ and $M$ is any integer divisible
by $2\lcm(m_{ij},m_{jk},m_{ik})^2$. Choose such positive $M$ common
for all triples $\{i,j,k\}$ and let $F_{ij}=f_{ij}^{M/m_{ij}}$. Then
$$
  \div (F_{ij}) = M(P_i) - M(P_j) \qquad (\{i,j\}\subset\{1,2,3,4\})
$$
and
\beq\label{prodF}
  \frac{F_{ij}}{F_{ij}(P_k)}
  \frac{F_{jk}}{F_{jk}(P_i)}
  \frac{F_{ki}}{F_{ki}(P_j)} = 1 \qquad (\{i,j,k\}\subset\{1,2,3,4\})\>.
\eeq
Also define elements of $\KC$ by
$$
  S_{i,j,k} = \Bigl\{ \frac{F_{ki}}{F_{ki}(P_j)}, \frac{F_{ij}}{F_{ij}(P_k)} \Bigr\}
            = \Bigl\{ \frac{F_{ij}}{F_{ij}(P_k)}, \frac{F_{jk}}{F_{jk}(P_i)} \Bigr\}
            = \Bigl\{ \frac{F_{jk}}{F_{jk}(P_i)}, \frac{F_{ki}}{F_{ki}(P_j)} \Bigr\} \>.
$$
The asserted relation~\eqref{Crelation}, but with $c_i$ replaced by
\begin{equation}
\label{Ciprimes}
c_i'=\frac{\text{lcm}(m_{lj},m_{lk},m_{jk})}{m_{lj}m_{lk}m_{jk}} \qquad(\{i,j,k,l\}=\{1,2,3,4\})
\,,
\end{equation}
then follows from the statement that, in $K_2^T(C)/$torsion,
\beq\label{saltsum}
  S_{1,2,3} - S_{1,2,4} + S_{1,3,4} - S_{2,3,4} = 0 \>.
\eeq
To prove (\ref{saltsum}), let $f\!=\!F_{14}/F_{14}(P_2)$, $g\!=\!F_{24}/F_{24}(P_3)$ and $h\!=\!F_{34}/F_{34}(P_1)$.
Rescaling if necessary, we may also assume that
$F_{12}=fg^{-1}, F_{23}=gh^{-1}$ and $F_{31}=hf^{-1}$. Finally, let
$$
  \alpha = f(P_3), \>\>
  \beta=   g(P_1), \>\>
  \gamma = h(P_2) \>.
$$
Now we apply (\ref{prodF}) for
$\{i,j,k\}\!=\!\{1,2,4\}$, $\{1,3,4\}$, $\{2,3,4\}$ and $\{1,2,3\}$. This
gives, respectively, that
$$
  (fg^{-1})(P_4) = \beta^{-1}, \quad
  (gh^{-1})(P_4) = \gamma^{-1}, \quad
  (hf^{-1})(P_4) = \alpha^{-1} \quad \text{and} \quad
  \alpha\beta\gamma = 1\>.
$$
Finally, using all of these we expand the symbols $S_{i,j,k}$ in terms
of the 8 generators,
$$
  v = (
    \{f,g \}, \{f,h\}, \{g,h\}, \{f,\beta\}, \{g,\alpha\},
    \{g,\beta\}, \{h,\alpha\}, \{\alpha,\beta\}
  ) \>.
$$
We find
$$
  \begin{array}{rll}
  S_{1,2,3}\> = & \!\!\!\{fg^{-1}/\alpha, gh^{-1}/\beta\}\!\!\! & = \>(1,-1,1,-1,1,1,-1,1)\cdot v \>, \cr
  -S_{1,2,4}\> = & \!\!\!\{f, g/\beta\} & = \>(-1,0,0,1,0,0,0,0)\cdot v \>, \cr
  S_{1,3,4}\> = & \!\!\!\{f/\alpha, h\} & = \>(0,1,0,0,0,0,1,0)\cdot v \>, \cr
  -S_{2,3,4}\> = & \!\!\!\{g, h/\gamma\} & = \>(0,0,-1,0,-1,-1,0,0)\cdot v \>. \cr
  \end{array}
$$
Therefore the left-hand side of (\ref{saltsum}) reduces to one
symbol $\{\alpha,\beta\}$. This symbol can be also rewritten in terms
of the original functions. One shows that
\begingroup
$$
  \def\Q#1#2#3#4{\Bigl(\frac{f_{#1#2}#3}{f_{#1#2}#4}\Bigr)^{\frac{M}{m_{#1#2}}}}
  \{\alpha,\beta\} = \pm \Bigl\{ \Q kl{(P_i)}{(P_j)}, \Q il{(P_j)}{(P_k)} \Bigr\} \>\>.
$$
\endgroup
for all combinations of four distinct indices $ i $, $ j $, $ k $
and $ l $ with $|i-k|=|j-l|=2$. In any case
this element comes from $K_2$ of a number field, which is torsion.
Hence the asserted relation holds.

It remains to show that the numbers $c_i$ and $c_i'$ defined by~\eqref{Cis} and~\eqref{Ciprimes}
are proportional. Equivalently, if $\c_i=v_p(c_i)$ and $\c_i'=v_p(c_i')$
denote the valuations of $c_i$ and $d_i'$ at some prime~$p$, then we must show
that $\c_i-\c_i'$ is independent of $i$.  The numbers $\c_i$ and $\c_i'$ are
given by
\begin{equation}
\label{471}
\c_i=\min_{j\ne i}\{\nu_{ij}\},\qquad\c_i'=\max_{j,k\ne i}\{\nu_{jk}\}
  \,-\sum_{j,k\ne i,\;j<k}\nu_{jk}\qquad(i=1,2,3,4)\,,
\end{equation}
where $\nu_{ij}$ denotes $v_p(m_{ij})$.  If we use the six numbers $\nu_{ij}$
to label the edges of a tetrahedron $T$ in the obvious way, then the relation~\eqref{L}
says that the two largest labels of the sides of any face of~$T$ are
equal.  From this we find easily that there are two possibilities: if three
incident edges of $T$ have the same label, then, possibly after renumbering,
we have $\nu_{12}=\nu_{13}=\nu_{14}=a$, $\nu_{23}=\nu_{24}=b$, $\nu_{34}=c$
for some integers $a\ge b\ge c$, while if this does not happen then, again
up to renumbering, we have $\nu_{13}=\nu_{14}=\nu_{23}=\nu_{24}=a$,
$\nu_{12}=b$, $\nu_{34}=c$ with $a > b,\,c$.  From~\eqref{471} we then find
that $(\c_1,\ldots,\c_4)$ and $(\c'_1,\ldots,\c'_4)$ are given by $(a,b,c,c)$
and $(-b-c,-a-c,-a-b,-a-b)$ in the first case and by $(b,b,c,c)$ and
$(-a-c,-a-c,-a-b,-a-b)$ in the second case, so that in both cases
$\c_i-\c_i'=a+b+c$ for all $i$.  This completes the proof of Proposition~\ref{tetrel}.
\end{proof}

\begin{corollary}\label{deppts}
Let $C$ be a curve defined over $\Q$ and 
$P_1,\ldots,P_n\in C(\Q)$ points such that all $(P_i)-(P_j)$ are torsion. 
Then the subspace of $K_2^T(C)\tensor\Q$ generated by all elements $\{P_i,P_j,P_k\}$
is already generated by those of the form $\{P_1,P_i,P_j\}$.
\end{corollary}

The corollary implies that the space spanned by the $n(n-1)(n-2)$ symbols $\{P_i,P_j,P_k\}$,
which already by part~(2) of Proposition~\ref{consprop} had dimension at most $\binom n3$, 
in fact has dimension at most $\binom{n-1}2$.

\begin{remark}\label{pointsqbar}
If $C/\Q$ is a curve, it is sometimes convenient to consider points $P_i$
in $C(\Qbar)$. Then we have to work in \KCQbar, which is defined as
the kernel of the tame symbol (given by \eqref{tame}),
$$
  \KCQbar = \Ker\bigg( K_2(\Qbar(C)) \to \bigoplus_{ x \in C(\Qbar) } \Qbar^* \bigg) \> .
$$
The product formula~\eqref{productformula} still holds for elements
in $ K_2(\Qbar(C)) $.
All the results in this section remain true in this context.
Moreover, since the inclusion of $\Q(C)$ into $\Qbar(C)$ induces a map from
$K_2(\Q(C))$ to $K_2(\Qbar(C))$ and
the tame symbol on both groups is given by the same formula,
we can check if an element in $K_2(\Q(C))$ lies in
$\KC$ using the tame symbol on its image in $K_2(\Qbar(C))$.

In fact, the map $ K_2(\Q(C)) \to K_2(\Qbar(C)) $ has torsion kernel
and we may identify $ K_2(\Q(C)) \tensor \Q $ with
the subspace of $ K_2(\Qbar(C)) \tensor \Q$ on which $ \Gal(\Qbar/\Q) $
acts trivially.  The same statements hold when we replace $ K_2(\Q(C)) $ by $ K_2^T(C) $
and $ K_2(\Qbar(C)) $ by $ K_2^T(C_\Qbar) $

\end{remark}

The final result of this section, a strengthening of Corollary~\ref{deppts}, says that we cannot construct any more elements 
of $K_2^T(C_\Qbar) \tensor \Q$ using only functions whose zeros and poles differ by torsion divisors
than those already given by Construction~\ref{constors}.

\begin{proposition}\label{Wproposition}
Let $ S \subseteq C(\Qbar)$ 
and $ P_0 \in S $ be such
that $ (P)-(P_0) $ is a torsion divisor for all $ P $ in $ S $. Let $V$ be the subspace of 
$K_2^T(C_{\Qbar})\tensor\Q$ generated by all elements $\{P,Q,P_0\}$ with $P,Q\in S$
and $W$ the subspace of $K_2(\Qbar(C))\tensor\Q$ generated by all symbols $\{f,g\}$ with
$\,\text{\rm div}(f)\,$ and $\,\text{\rm div}(g)\,$ supported in $S$. Then
$$ W\,\cap\,K_2^T(C_{\Qbar})\tensor\Q\;=\;V\,.$$   
\end{proposition}

\begin{proof}
Let $\xi=\sum_i\{f_i,g_i\}$ be an element of $W\cap\,K_2^T(C_{\Qbar})\tensor\Q$.
Replacing $ S $ with the finite set of all zeros and poles of the $f_i$ and $g_i$,
by assumption there exist an integer $N>0$ and, for each $P\in S\setminus\{P_0\}$,
a function $f_P\in\Qbar(C)^*$
with $\,\text{div}(f_P)=N(P)-N(P_0)$. Then each $f_i^N$ or $g_i^N$
is, up to a scalar factor, a multiplicative combination of the functions $f_P$,
so $N^2\xi$ is a linear combination of symbols of the form $\{f_P,f_Q\}$, $\{f_P,c\}$
and $\{c,c'\}$ with
$c,\,c'\in\Qbar^*$.  But $\{c,c'\}=0$ 
in $K_2^T(\Qbar(C))\tensor\Q$ and $\{f_P,f_Q\}$ is the sum of $\{f_P,f_Q(P)\}+\{f_P(Q),f_Q\}$
and a multiple of $\{P,Q,P_0\}\in V$, so $\xi$ can be written as 
$\sum_P \{f_P,c_P\}+\xi'$
with $\xi'\in V$ for some constants $c_P\in\Qbar^*$.
Since the tame symbol of $\{f_P,c_P\}$ at $x$ is $c_P^{-N}$ for $x=P$ and trivial for
all other $x\in C(\Qbar)\setminus\{P_0\}$, the fact that some multiple of $\xi$
(and hence also of $\xi-\xi'$) has trivial tame symbol everywhere implies that each
$c_P$ is a root of unity. Hence $\xi\in V$. 
\end{proof}


\section {Torsion divisors on hyperelliptic curves} 

Now suppose that $C$ is a hyperelliptic curve of genus $g\ge 1$ defined over $\Q$.
The hyperelliptic involution on $C$ determines
a double cover $\phi: C\to \P^1$ ramified at $2g+2$ points, the fixed
points of the involution. Assume that one of these points $\infty\iC$ is defined
over \Q. After a change of coordinates on $\P^1$ we can assume that
$\phi(\infty\iC)=\infty\si{\P^1}$. Then $C$ gets a model
\beq\label{hyperelleq}
  y^2 = c_{2g+1} x^{2g+1} + c_{2g} x^{2g} + \ldots + c_1 x + c_0, \qquad c_{2g+1}\ne 0,
\eeq
where the polynomial on the right has coefficients in \Q and has no
multiple roots. This equation can be seen either as defining a double
cover of $\P^1$ or as a curve (singular for $g>1$) in $\P^2$ whose
normalization is $C$. The point $\infty\!=\!\infty\iC$ is then the unique point at infinity
of this normalization. The cover $\phi$ is given by the function $x$ on
$C$ and its ramification points are
$\infty$ and the $T_\a\!=\!(\alpha,0)$, where $\alpha$
runs through the roots of the right-hand side of (\ref{hyperelleq}) in $\Qbar$.

We will look for points $P$ such that the divisor $(P)-(\infty)$ is
$m$-torsion for some~$m$.
Such points will be called ($m$-){\em torsion points}. Of course, if $P, Q$ are torsion
points, then $(P)-(Q)$ is also a torsion divisor.
If we succeed in constructing many torsion points $P_i$, then we get
many elements $\{\infty, P_i, P_j\}\in\KC/\torsion$ by using
Proposition~\ref{consprop}.
Here are some examples of curves with torsion points:

\ExampleNamed{2-torsion}\label{2torsex}
With notation as before,
any non-trivial difference $(P)-(Q)$ for $P$ and $Q$  among  $\infty$ and the $T_\a$ is 2-torsion:
$$
  2(T_\a)-2(T_\b)=\div\Bigl(\frac{x\!-\!\a}{x\!-\!\b}\Bigr), \qquad
  2(T_\a)-2(\infty)=\div(x\!-\!\a)\>.
$$
So the $T_\a$ are 2-torsion points, though not necessarily defined over \Q.

Conversely, suppose given an arbitrary curve $C$ over $ \Q $
and two distinct points
$P,Q\in C(\Q)$ such that the divisor $(P)-(Q)$ is 2-torsion.
Say $\div(\phi)=2(P)-2(Q)$ with $\phi$ defined over $\Q$.
Then $\phi: C\to\P^1$ is a double cover defined over $\Q$.
It follows that $C$ is hyperelliptic and admits a model
(\ref{hyperelleq}) with $Q=\infty$, $P$ among the $T_\a$.

Unfortunately, we cannot use 2-torsion points alone to obtain interesting elements in $\KCQbar/\torsion$ by the construction of Section~\ref{sTorssymb}:
the calculation 
\beq*%
\begin{array}{r}
  \displaystyle
    \Bigl\{ \frac{x-\b}{\a-\b} , \frac{\b-\a}{x-\a}\Bigr\}
  = \left\{ \frac{x-\a}{\b-\a} , \frac{x-\b}{\a-\b} \right\}
  = \left\{ \frac{x-\a}{\b-\a} , 1-\frac{x-\a}{\b-\a} \right\}
  = 0
\end{array}
\eeq*%
shows that  $\{\infty,T_\a,T_\b\}$  is trivial 
for any $\a$ and $\b$, and 
Proposition~\ref{tetrel} (or rather, its analogue over $\Qbar$) then shows that the 
elements  $\{T_\a,T_\b,T_\c\}$ also vanish.  The 2-torsion points can nevertheless
be used, but only in combination with the torsion points of other orders which
we describe next.

\ExampleNamed{($\mathbf{2g\!+\!1}$)-torsion}
\label{ex2g1tors}
Assume that the hyperelliptic equation (\ref{hyperelleq})
has the special form
\beq\label{y2_2g1}
  y^2 = c\> x^{ 2g+1 } + (b_g x^g + \ldots + b_1 x + b_0)^2 
\eeq
for some $ c \neq 0 $.  We can scale $x$ and $y$ to make $c=-1$ if desired,
which allows a uniform treatment together with Example~\ref{ex2g2tors} below.
The substitution $y\!\mapsto\!y+\sum b_ix^i$ transforms
this equation to
\beq\label{y2y_2g1}
  y^2 + 2(b_g x^g + \ldots + b_1 x + b_0)y = c x^{2g+1} \>.
\eeq
The point $O=(0,0)$ (corresponding to $(0,b_0)$ on the curve (\ref{y2_2g1}))
lies on this curve and the function $y$
has divisor $(2g+1)(O)\!-\!(2g+1)(\infty)$, so $O$
is a $(2g+1)$-torsion point.

In fact, any curve $C$ as in (\ref{hyperelleq}) with a non-trivial rational
$(2g+1)$-torsion point $O$ is isomorphic to a curve of the form (\ref{y2_2g1})
and hence (\ref{y2y_2g1}).
Indeed, given such a curve $C$, consider the local system
$$
L = H^0(C, \cO\iC((2g+1)\infty)).
$$
By Riemann-Roch, the dimension of this as $\Q$-vector space
is $ 1 -g + 2g+1  = g + 2 $. It is easy to see that
$\{1,x,\ldots,x^g,y\}$ is a basis.
So for a point $ O=(x\iO,y\iO)\in C(\Q)$ the divisor
$ (2g+1)(O)-(2g+1)(\infty)$ is principal
if and only if there is an element of $ L  $,
$$
  h(x,y) =  y - \tilde{b}_g x^g - \ldots - \tilde{b}_1 x - \tilde{b}_0\>,
$$
which vanishes to order exactly $ 2g+1 $ at $ O $.
(Note that necessarily it cannot be expressed in terms of $ x $ alone
because the order of such a function is at most $ 2g $.)
Then the curve defined by $h(x,y)=0$ can have only the point $(x\iO,y\iO)$
in common with $C$, and after substituting
$y = \tilde{b}_g x^g + \ldots + \tilde{b}_1 x + \tilde{b}_0$ in~(\ref{hyperelleq}),
we see that the equation of the curve is of the form
$$
  y^2 = c (x-x\iO)^{ 2g+1 } + (\tilde{b}_g x^g + \ldots + \tilde{b}_1 x + \tilde{b}_0)^2 \>.
$$
After a translation we can assume $x\iO=0$
and the equation of the curve becomes of the form (\ref{y2_2g1})
with $O=(0,b_0)$.

\ExampleNamed{($\mathbf{2g\!+\!2}$)-torsion}
\label{ex2g2tors}
Similar to the previous example, a curve
$$
  y^2 = -c^2 x^{ 2g+2 } + (c x^{g+1} + b_g x^g + \ldots + b_1 x + b_0)^2 \>, 
$$
with $b_g, c\ne 0$ is of the form $(\ref{hyperelleq})$
and it has a model
$$
  y^2 + 2(c x^{g+1} + b_g x^g + \ldots + b_1 x + b_0)y = -c^2 x^{ 2g+2 } \>.
$$
In this model  $\div(y) = (2g+2)(O)-(2g+2)(\infty)$, so the curve
has a $(2g+2)$-torsion point $O\!=\!(0,0)$.
Scaling $ y $ we can assume that $ c=1 $, hence $-c^2=-1$, allowing
a uniform treatment together with Example~\ref{ex2g1tors}.

As in Example~\ref{ex2g1tors}, one can use the linear system
$$
  L = H^0(C, \cO\iC((2g+2)\infty))
$$
to show that any curve (\ref{hyperelleq}) with a rational $(2g+2)$-torsion
point has a model of this form.


\section {Elements of $K_2$ for hyperelliptic curves} 

\begingroup
\def\sym#1#2{\Bigl\{\>#1\>,\>#2\>\Bigr\}}
\def\symf#1#2#3#4{\sym{\frac{#1}{#2}}{\frac{#3}{#4}}}
\def\S#1,#2;{\displaystyle\Bigl\{\,#1,\,#2\,\Bigr\}}%

In this section we study elements of $\KC$ constructed from torsion
points on a non-singular hyperelliptic curve $C/\Q$ of genus
$ g $. We show that using only
2-torsion points is not sufficient to construct non-torsion elements.  So we will use
curves which also have rational torsion points of order $2g+1$ or $2g+2$.

As discussed in Examples \ref{ex2g1tors} and \ref{ex2g2tors},
such a curve $C/\Q$ can always be given by an equation of the form
\beq\label{eqdeff}
  y^2 + f(x)\,y + x^d = 0 
\eeq
in one of the following two cases:
\beq\label{eqdeff2}
\begin{array}{llll}
\qquad\text{(a)} & d=2g+1,  & f(x)= b_g x^g + \ldots + b_1 x + b_0\>; \cr
\qquad\text{(b)} & d=2g+2,  & f(x)= 2x^{g+1}+b_g x^g + \ldots + b_1 x + b_0\text{ and } b_g \neq 0\>. \cr
\end{array}
\eeq%
This curve is isomorphic to $y^2 = t(x)$,
where $t(x)$ is the 2-torsion polynomial,
$$
  t(x) = -x^d+{f(x)^2}/4 \>.
$$
Since we want $C$ to be non-singular, we assume that $t(x)$ has no
multiple roots, and, in particular, that $b_0\ne 0$.

On the curve \eqref{eqdeff} we have a rational point $O=(0,0)$ and the divisor
$d(O)\!-\!d(\infty)$ is principal, namely equal to $\div (y)$.
There is also a ``reflected'' $ d $-torsion point
$$
  O' = (0,-f(0)),
$$
the image of $O$ under the hyperelliptic involution
$(x,y)\mapsto(x,-y-f(x))$.

Moreover, every rational root $\alpha$ of $t(x)$
gives a rational 2-torsion point $T_\alpha=(\alpha,-f(\alpha)/2)$.
All the pairwise differences of the points $O$, $O'$, $\infty$ and the $T_\alpha$
are torsion divisors. So the results of Section \ref{sTorssymb} apply and we
get elements of $\KCtors$, such as
\beq*
  \{\infty, O,T_\alpha\}, \quad
  \{\infty,O, O'\}, \quad
  \{T_\alpha, T_\beta, T_\gamma\}, \quad\ldots\quad.
\eeq*%
If we replace  $\KC$ with $\KCQbar$ then, by
Remark~\ref{pointsqbar}, we can also use points $T_\alpha$ for
any root $\alpha\in\Qbar$ of the 2-torsion polynomial $t(x)$, not necessarily
rational.
The following proposition gives relations in $ K_2^T(C_\Qbar) \tensor \Q $ among these elements.
We mention that the two relations involving $ \{\infty,O,O'\} $ were first observed
in computer calculations at the regulator level. 

\begin{proposition}\label{Vproposition}
Assume $ g \geq 2 $, and let $V$ be the \Q-subspace of $ \KCQB\tensor\Q $ generated by the elements
of the form $\{P_i, P_j, P_k\}$ where $P_i, P_j, P_k$ run through the points
$\infty, O, O'$ and the $T_\alpha$ for $\alpha\in\Qbar$ a root of $t(x)$.
Then $V$ is already generated by the elements of the form $\{\infty, O, T_\alpha\}$.
More precisely,
if we write $ d = 2g+\varepsilon $  with $ \varepsilon =1$  or 2, then we have the relations
\begin{align*}
&\qquad\qquad \{\infty, T_\a, T_\b\} = 0 \> , \cr
&\qquad\qquad \{\infty, O, T_\a\} + \{\infty, O', T_\a\} = 0  \cr
\intertext{and}
&\qquad\qquad \e\{\infty, O, O'\} = \sum_\a \{\infty, O, T_\a\} \>,
\end{align*}
together with the symmetry properties of Proposition~\ref{consprop} and the tetrahedron relation of Proposition~\ref{tetrel}.

Furthermore, if $ f(0)^2=1 $, then
\[ \{-f(0)y,-x\} \]
is in $  K_2^T(C) $, and its class $ \specelt $ in $ K_2^T(C)/\torsion $ satisfies 
\[
d\,\specelt = \e\{\infty, O, O'\}
\>.
\]
\end{proposition}

\begin{proof}
The proof of the first relation was already given in Example~\ref{2torsex}.
The other three relations involve
$ \{\infty, O, T_\alpha\} $,
$ \{\infty, O', T_\alpha\} $ and
$ \{\infty, O, O'\} $.
To make these elements
explicit from the definitions in Construction~\ref{constors} and Proposition~\ref{consprop}
we need to determine the
order $ d' $ of $ (O)-(O') = 2(O)-2(\infty) $ and
the order $ d_\a'' $ of $ (O)-(T_\a) = (T_\a) -(O') $
in $ \Pic^0(C) $.

In fact, we have
\begin{equation*}
 d' = d/\e  \qquad  \text{and} \qquad  d_\a'' = 2 d / \e  \> .
\end{equation*}
Namely, if $ d=2g+1 $ then $ d'= d $ and if $ d=2g+2 $ then $ d'=d/2 $.
And if $ d=2g+1 $, then $ (O) - (T_\a) = [(O)-(\infty)] - [(T_\a)-(\infty)] $
must have order $ d_\a'' = 2d $, but if $ d=2g+2 $ the order
can be $ d $ or $ d/2 $.  However, it can be $ d/2 = g+1 $
only if $ (g+1)[(O)-(\infty)] = (T_\a)-(\infty) = (\infty) - (T_\a) $, so
that $(g+1)(O) + (T_\a) - (g+2)(\infty) $ is the divisor of some
non-zero $ h $ in $ H^0(C, \cO\iC((g+2)\infty)) $.
Since $ g\geq 2 $, $ g+2 \leq 2g $, so $ h $ lies in $ H^0(C, \cO\iC(2g\infty)) $,
which has basis $ \{1,x,\dots,x^g\} $ (cf.~Example~\ref{ex2g1tors}).
But that implies that $ \div(h) $ is invariant under the hyperelliptic involution
whereas $(g+1)(O) + (T_\a) - (g+2)(\infty)$ is not because $ O \neq O' $.
So $ d_\a''= d/2 $ is impossible.

Using those values of $ d' $ and $ d_\a'' $ in the definitions we now see that
\begin{align}
\{ \infty, O, T_\a \}  & = \S \frac{y}{-f(\alpha)/2} , \frac{x-\alpha}{-\alpha}; \>,\label{NEWlabel}\\
 \{ \infty, O', T_\a \} & = \S \frac{-f(x)-y}{-f(\alpha)/2}, \frac{x-\alpha}{-\alpha}; \nonumber\\
\intertext{and}
 \e \{ \infty , O , O' \} & = \S \frac{y}{-f(0)}, \frac{y+f(x)}{f(0)}; \> ,\nonumber
\end{align}
where we write elements in $ \KCQbar $ for their classes in $ \KCQbar/\torsion $.

The second relation in the proposition now follows from 
\beq*
\Bigl\{
\frac{y}{-f(\alpha)/2} , \frac{x-\alpha}{-\alpha}
\Bigr\}
+
\Bigl\{
\frac{-f(x)-y}{-f(\alpha)/2} , \frac{x-\alpha}{-\alpha}
\Bigr\}
=
\Bigl\{
\frac{y(f(x)+y)}{-f(\alpha)^2/4} , 1 - \frac{x}{\alpha}
\Bigr\}
\cr
=
\Bigl\{
\frac{- x^d}{- \alpha^d} , 1 - \frac{x}{\alpha}
\Bigr\}
=
d
\Bigl\{
\frac{x}{\alpha} , 1 - \frac{x}{\alpha}
\Bigr\} = 0 \>.
\eeq*%
In order to prove the third relation, we notice that
\beq\label{eltmul2}
\begin{array}{r}
  \displaystyle
  2\symf{y}{-f(\a)/2}{x-\a}{-\a} =
  \sym{\frac{y^2}{f(\a)^2/4}}{1-\frac x\a} =
  \sym{\frac{-y^2}{-\a^d}}{1-\frac x\a}
\cr\displaystyle
  = \sym{\frac{y^2}{\a^d}}{1-\frac x\a} - d \sym{\frac{x}\a}{1-\frac x\a} =
  \sym{\frac{y^2}{x^d}}{1-\frac x\a} \> ,
\end{array}
\eeq
so that taking the sum over all roots $\alpha\in\Qbar$ of $t(x)$
gives
\begin{equation}
\label{oopelement}
\begin{array}{ll}
 &   \displaystyle\sum_{t(\alpha)=0} \>\> \S  \frac{y^2}{x^d}, 1-\frac x{\alpha};
   = \S \frac{y^2}{x^d},  \prod_{t(\alpha)=0}\Bigl(1-\frac x{\alpha}\Bigr);
 \cr
   &\qquad\qquad
     =  \S \frac{y^2}{x^d}, \frac{t(x)}{t(0)};
     =  \S  \frac{y^2}{x^d}, \frac{(y+f(x)/2)^2}{f(0)^2/4};      
 \>,
\end{array}
\end{equation}
where we used that $t(x)=(y+f(x)/2)^2$.
Also, already in $ \KC $,
\begin{equation}
\label{lastidentity}
\begin{array}{rl}
2  \S \frac{y}{-f(0)}, \frac{y+f(x)}{f(0)}; \!\!\!\!
& =
2  \S \frac{y}{-f(0)}, \frac{y\,(y\!+\!f(x))}{f(0)^2};
\cr
 & =
  \S \frac{y^2}{f(0)^2}, \frac{-x^d}{f(0)^2};   = 
  \S  \frac{y^2}{x^d}, \frac{-x^d}{f(0)^2};
\> ,
\end{array}
\end{equation}
so that $2 \e \{\infty,O,O'\} - 2 \sum_\a \{\infty, O, T_\a \} $ is equal to the class in $ \KCQbar/\torsion $
of
\beq
\label{reallylastidentity}
\qquad
  \S  \frac{y^2}{x^d}, \frac{-x^d}{f(0)^2};  -
  \S  \frac{y^2}{x^d}, \frac{(y+f(x)/2)^2}{f(0)^2/4};        =
  \S  \frac{y^2}{x^d}, \frac{-x^d}{(2y+f(x))^2};   
\,.
\eeq%
That this last element is zero in $ \KCQbar $ can be seen by
applying the identity
\begin{align*}
 \S \frac{a}{1-a}, \frac{a(a-1)}{(1-2a)^2};
  & =
 \S \left(\frac{1-a}{a}\right)^2, 1-\left(\frac{1-a}{a}\right)^2;
  -
 \S \frac{a}{1-a}, -\frac{a}{1-a};
\cr
  & \qquad\qquad\qquad\qquad\qquad\qquad -2 \S a, -a; +2 \S 1-a, a;
  = 0
\>,
\end{align*}
valid in $ K_2 $ of any field if $ a $, $ 1-a $ and $ 1-2a $
are non-zero, to $ a=-y/f(x) $, so that $ 1-a = -x^d/f(x)y $
and $ 1-2a = (2y+f(x))/f(x) $.  This comples the proof of the
third relation.

%

Finally, if $ f(0)^2 =1 $, then all the tame symbols $T_P$ of
$ \{-f(0)y,-x\} $ are trivial
for $ P $ different from $ O=(0,0) $, $ O'=(0,-f(0)) $ and $ \infty $,
and for those three points we find that
\beq*
\begin{array}{lll}
T_{(0,0)}(\{-f(0)y,-x\}) & = & (-1)^d \dfrac{-f(0)y}{(-x)^d} \Bigm |_{(0,0)} = 1\>,
\cr
T_{(0,-f(0))} (\{-f(0)y,-x\}) & = & (-1)^0 \dfrac{-f(0)y}{(-x)^0} \Bigm |_{(0,-f(0))} = f(0)^2 = 1\>,
\cr
T_\infty(\{-f(0)y,-x\}) & = & 1 \>.
\end{array}
\eeq*%
The first equation uses that $ y(y+f(x)) = -x^d $, so that
the function $ -y/x^d $ equals $ 1/(y+f(x)) $ and therefore
assumes the value $ 1/f(0) $
at $ (0,0) $.  The third equation follows from the first two and the
product formula.  
Then~\eqref{lastidentity} shows that
\[
 2 d \{-f(0)y,-x\} = \{ y^2, - x^d \} = 2 \S \frac{y}{-f(0)}, \frac{y\!+\!f(x)}{f(0)};
\]
in $ \KC $, so that $ 2 d \specelt = 2 \e \{\infty, O, O'\} $
in $ \KC/\torsion $.  This finishes the proof.
\end{proof}

\begin{remark}\label{genusone}
When $ g=1 $, the statements of Proposition~\ref{Vproposition} still hold,
but when $ d=2g+2=4$, the identity $ \e \{\infty, O, O' \} = \sum_\a \{\infty, O, T_\a\} $
is replaced by
\[ 2\{\infty, O, O'\} = \{\infty, O, T_\a\} +  \{\infty, O, T_{\a'}\} + 2 \{\infty, O, T_\b\} \>, \]
where $ T_\a $, $ T_{\a'} $ and $ T_\b $ are the points of order~2, and
$ 2(O)-2(\infty) = (T_\b)-(\infty) $ in $ \Pic^0(C_\Qbar) $.
The proof is the same, using that
$ 2 \{\infty,O,O'\} =  \{ \frac{y}{-f(0)}, \frac{y+f(x)}{f(0)}\} $,
$ \{ \infty, O, T_\c \}   = \{ \frac{y}{-f(\c)/2} , \frac{x-\c}{-\c} \} $
for  $ \c = \a $ or $ \a' $ but 
$ 2 \{ \infty, O, T_\b \}   = \{ \frac{y}{-f(\b)/2} , \frac{x-\b}{-\b}\} $,
and similarly for $ O' $ instead of $ O $ in the last two identities.
(Here again we write elements in $ \KCQbar $ for their classes in $ \KCQbar/\torsion $.)
In~\eqref{reallylastidentity} one takes 
$ 
  4 \{\infty,O,O'\} - 2\{\infty,O,T_\a\} - 2\{\infty,O,T_{\a'}\} - 4 \{\infty,O,T_\b\}
$ 
as the left-hand side, whereas the right-hand side vanishes as before.
\end{remark}

Now we return from $ K_2^T(C_\Qbar)/\torsion$ to \KC/\torsion, the group that we are interested in.
Here we have elements $\{\infty, O, T_\a\}$, where $\alpha$ is a rational
root of $t(x)$ but, in fact, we have more.
As made explicit in \eqref{NEWlabel}, $\{\infty, O, T_\alpha\}$ is
the class of an element of $ \KC $, namely 
$ \{ \frac{y}{-f(\alpha)/2}, \frac{x-\alpha}{-\alpha} \} $.
(If $ g=1 $ and $ d=4 $ this statement has to be slightly modified
according to Remark~\ref{genusone}.)
If $m(x)$ is a rational factor of the 2-torsion polynomial $t(x)$, then
the element 
$
   \sum_{m(\alpha)=0}
   \{ \frac{y}{-f(\alpha)/2}, \frac{x-\alpha}{-\alpha} \}
$
(where the sum is taken over the roots of $m(x)$ in $\Qbar$)
in $ \KCQbar $
can be shown to come from $\KC$.  But if we multiply it by 2 and
use~\eqref{eltmul2}, then, by a calculation similar to~\eqref{oopelement},
we get explicitly that
\beq
\label{normcalculation}
2 \sum_{m(\alpha)=0} \>\> \Bigl\{
      \frac{y}{-f(\a)/2},1\!-\!\frac x{\alpha}
    \Bigr\} =
   \Bigl\{
      \frac{y^2}{x^d},\frac{m(x)}{m(0)}
    \Bigr\} \>.
\eeq%
This computation shows that we get an element of $\KC$ for every
rational irreducible factor $m(x)$ of the 2-torsion polynomial $t(x)$,
not just for the linear ones.
Thus we get $k$ explicit elements of $\KC/\torsion$,
where $k$ is the number of irreducible rational factors of the 2-torsion polynomial
of the curve in~\eqref{eqdeff}.  This gives a map from $ \Z^k $ to $ \KC/\torsion $.
In summary, we have:

\Construction\label{maincons}
Let $C/\Q$ be given by (\ref{eqdeff}), (\ref{eqdeff2}).
Let $m_1,\dots,m_k$ be the
irreducible factors in $\Q[x]$ (up to multiplication by $ \Q^* $)
of the 2-torsion polynomial $t(x) = -x^d+f(x)^2/4$.
To each of them we associate an element of $\KC/\torsion$,
\beq\label{Mjdef}
  M_j = \text{the class of } \Bigl\{
      \frac{y^2}{x^d},\frac{m_j(x)}{m_j(0)}
    \Bigr\},\qquad 1\le j\le k \>.
\eeq
Extending by linearity gives a map
\beq\label{LinEmbK2}
  \begin{array}{ccl}
     \Z^k & \longrightarrow & \quad \KC/\torsion \cr
  (n_1,\dots,n_k) & \longmapsto & \displaystyle \sum_{j=1}^k n_j M_j =  \text{the class of } \Bigl\{\frac{y^2}{x^d},\prod_{j=1}^k\frac{m_j(x)^{n_j}}{m_j(0)^{n_j}}\Bigr\} \>.
  \end{array}
\eeq

In two special cases, the statement of Construction~\ref{maincons} can be refined a little.  One case is
when $ f(0)^2 = 1 $, when we can use the class $ \specelt $ of Proposition~\ref{Vproposition} to replace
the $ \Z^k $ in~\eqref{LinEmbK2} by a larger lattice of the same rank.  The other situation occurs only
if $ d=2g+2 $, when it turns out that there is always a relation among the $ M_j $
(also originally discovered during calculations of the regulators)
and hence the rank of the image is at most $ k-1 $.

\begin{proposition}\label{relprop}
Let all notation be as in Construction~\ref{maincons}.  
\begin{enumerate}
\item\label{specelt}
If $ f(0) = \pm 1 $, then the class $\specelt$ in $ \KC/\torsion $ of $ \{-f(0)y,-x\}$ satisfies
$$
  2d\specelt = \sum_{j=1}^k M_j \>.
$$
\item\label{urelation}
If $d=2g+2$, so that $ 4t(x) = (f(x) - 2x^{g+1})(f(x) + 2x^{g+1}) $,
and
$$
     f(x) - 2x^{g+1} = m_1\dots m_l \>,
\qquad
     f(x) + 2x^{g+1} = m_{l+1}\dots m_k
$$
with $ m_j $ irreducible, then the corresponding classes $ M_j $ satisfy
$$
  M_1 + \ldots + M_{l} = M_{l+1} + \ldots + M_k \>.
$$
\end{enumerate}
\end{proposition}

%

\begin{proof}
Relation~\eqref{specelt} follows from the two identities involving $ \e \{\infty,O,O'\} $ in Proposition~\ref{Vproposition},
\eqref{NEWlabel} and~\eqref{normcalculation}
(or the corresponding statements in Remark~\ref{genusone} if $ g=1 $ and $ d=4 $).

In order to prove~\eqref{urelation}, let $H=-y/x^{g+1}$, so that $H+1/H = f(x)/x^{g+1} $. Then the relation
follows immediately from the fact that, in $ \KC $,
\beq*
 \Bigl\{ \frac{y^2}{x^d}, \frac{f(x)-2x^{g+1}}{f(x)+2x^{g+1}} \Bigl\}  \,
=\!
 \Bigl\{ H^2, \frac{H-2+1/H}{H+2+1/H} \Bigl\}
 \, =  2 \Bigl\{ H, \frac{(1-H)^2}{(1+H)^2} \Bigl\}
\cr
\vphantom\int =  4 \{ H, 1-H \} - 4 \{ -H, 1-(-H) \} = 0 \>.
\eeq*%
\end{proof}

\begin{remark}
If we identify
$ K_2(\Q(C)) \tensor \Q $ with the $ \Gal(\Qbar/\Q) $-invariants
of $ K_2(\Qbar(C)) \tensor \Q$ as in Remark~\ref{pointsqbar},
and let $ W \subseteq  K_2(\Qbar(C)) \tensor \Q  $ 
be as in Proposition~\ref{Wproposition}
(with $ S = \{\infty,O, O', T_{\a_1},\dots,T_{\a_{2g+1}} \} $,
where $ \a_1,\dots,\a_{2g+1} $ are the roots of $ t(x) $ in $ \Qbar $,
and $ \infty $ playing the role of $ P_0 $), then we have that
\[
  W \cap K_2^T(C) \tensor \Q  = \langle M_1, \dots, M_k \rangle_\Q 
  \> ,
\]
where the $ M_j $ are as in Construction~\ref{maincons}, and
$ \langle \, \cdots \rangle_\Q $ indicates the
$ \Q $-vector space they generate.
Therefore the $ M_j $'s give us everything in 
$ K_2^T(C) \tensor \Q $ that we can get from combinations of
symbols $ \{ f, g \} $ where $ f $ and $ g $ are in $ \Qbar(C)^* $
and $ \div(f) $ and $ \div(g) $ are
supported in $ \{\infty, O, O', T_{\a_1},\dots,T_{\a_{2g+1}} \} $.

Namely, Proposition~\ref{Wproposition} shows that $ V = \KCQbar \tensor \Q \cap W $
is spanned by the $ \{ \infty, P, Q \} $ where $ P $
and $ Q $ are in $ \{O, O', T_{\a_1},\dots,T_{\a_{2g+1}} \} $,
and by Proposition~\ref{Vproposition} and its proof those elements
can be expressed in the $\{\infty, O, T_{\a_i} \}$
or the $ \{ \frac{y^2}{x^d},\frac{x-\a_i}{-\a_i} \} $.
Therefore, if $ k = \Q(\a_1,\dots,\a_{2g+1}) $,
then the action of $ \Gal(\Qbar/\Q) $ on $ V $ factors through $ \Gal(k/\Q) $,
and $ V^{\Gal(\Qbar/\Q)} $ is generated by the elements
\[
\sum_{\sigma \in \Gal(k/\Q)} \sigma \left(  \Bigl\{ \frac{y^2}{x^d},\frac{x-\a_i}{-\a_i}  \Bigr\}\right)
=
\sum_{\sigma \in \Gal(k/\Q)}   \Bigl\{ \frac{y^2}{x^d},\frac{x-\sigma(\a_i)}{-\sigma(\a_i)}  \Bigr\}
=  \delta_i \Bigl\{ \frac{y^2}{x^d}, \frac{n_i(x)}{n_i(0)}  \Bigr\}
\> ,
\]
where $ n_i(x) $ is the minimal polynomial of $ \a_i $ over $ \Q $
and $ \delta_i \cdot \deg(n_i(x)) = [k:\Q] $.
Because $ n_i(x) $ is an irreducible rational factor of $ t(x) $
this equals $ \delta_i M_j $ for some $ j $.
\end{remark}

\endgroup 


\section {Constructing good polynomials} 
\label{sFam}

Our goal is to use Construction \ref{maincons} to produce explicit families
of hyperelliptic curves $ C $ with as many elements of $\KC$ as possible.
This comes down to constructing polynomials of the right form,
which we will refer to as ``good'' polynomials. This is addressed in this section.

\begin{problem}
Construct (families of) polynomials of the form
\[
  t(x) = - x^d + f(x)^2/4
\]
which have many rational factors and no multiple roots.
Here $d\ge 5$, and we want that
$t(x)$ has degree $d-1$ for $ d $ even, and degree $ d $ for $ d $ odd.
\end{problem}

\Remark To get potentially interesting examples for the Beilinson
conjecture, we need at least $ g\!=\,\, $genus$(C)$ linearly independent
elements of $\KC/\torsion$. Thus
we want $t(x)$ to have at least $g$ rational factors if $d=2g+1$ is odd
and at least $g+1$ of them
if $d=2g+2$ is even (cf.~\eqref{urelation} of Proposition~\ref{relprop}).
Such a $ t(x) $ is what we will call a good polynomial.

We have seen in Examples~\ref{ex2g1tors} and~\ref{ex2g2tors} that,
for our purposes, we might just as well consider, for $ c \neq 0 $,
\[
t(x) =  cx^{2g+1} + f(x)^2 \text{ where } f(x) =  b_g x^g + \ldots + b_1 x + b_0
\]
when $ d=2g+1 $, and
\[
t(x)= -c^2x^{2g+2} + f(x)^2 \text{ where } f(x) = cx^{g+1}+b_g x^g + \ldots + b_1 x + b_0
\text{ with } b_g \neq 0 
\]
when $ d=2g+2 $.  It will be more convenient to use those non-normalized
versions because we can sometimes let $ c $ play a role in the construction
of such $ t(x) $.

For the remainder of this section
we keep the notations $ c $, $d$, $t(x)$ and $f(x)$ as above.
For $d=5$ and $d=6$ we will explain how to describe all
such $ t(x) $ that factor completely over the rationals.
For larger degrees we will give sporadic examples that factor
completely or nearly completely, and also infinite families of
good polynomials (with $ d = 2g+2 $) for all genera $ g $.

One way to produce good polynomials is to start with
a general polynomial $t(x)$ and force it to
have given rational roots. We illustrate this with one example:

\ExampleNamed{$\mathbf {d=5}$}
\label{ex5_2}
Take
$$
  t(x) = x^5 + (b_2 x^2 + b_1 x + b_0)^2 \>.
$$
and force the quintic to have two rational roots, so that we get 
(at least) 3 rational factors.

If $\alpha$ is a root of $t(x)$, then $-\alpha^5$ is a square, so
$\alpha=-k^2$ for some $k\in\Q$. We want $t(x)$ to have two distinct
rational roots $-m^2$ and $-n^2$, so
$$
-m^{10} + (b_2 m^4 - b_1 m^2 + b_0)^2 = 0
\quad\text{and}\quad
-n^{10} + (b_2 n^4 - b_1 n^2 + b_0)^2 = 0 \>.
$$
Taking roots yields
$$
m^5 = b_2 m^4 - b_1 m^2 + b_0
\quad\text{and}\quad
n^5 = b_2 n^4 - b_1 n^2 + b_0 \>.
$$
(We can take the positive signs by replacing $ m $ by $ -m $ or $ n $ by
$-n$ if necessary.)
Solve this linear system for $b_0$ and $b_1$ and
set $b_2=(k+m^2+mn+n^2)/(m+n)$. Finally, multiply by $(m+n)^2$ to
obtain
$$
  \def\sp{\!+\!}
  (m+n)^2 t(x) = (m \sp n)^2 x^5 +
    \Bigl(
      (k \sp m^2 \sp m n \sp n^2) x^2 + (k m^2 \sp m^2 n^2 \sp k n^2) x + k m^2 n^2
    \Bigr)^2 \>.
$$
This is a family on 3 parameters $m,n,k\in\Q$ although it gives only
a 2-dimensional family of curves $y^2=t(x)$: this equation is homogeneous
of multi-degree $(1,1,2,2,5)$ in $(m,n,k,x,y)$, so letting
$$
  m \mapsto \lambda m, \quad   n \mapsto \lambda n, \quad   k \mapsto \lambda^2 k\>,
$$
gives a polynomial which corresponds to an isomorphic curve.
Thus one can, for instance, assume that $m,n,k$ are integers or instead that,
say, $m=1$.
Note that $kmn\ne0$ (for otherwise the quintic has a double root $\alpha=0$)
and that $m\ne\pm n$.

\Remark\label{ex5_3}
We could also have forced $t(x)$ to have a third rational root
$\alpha=-l^2$. This condition gives a third linear equation for
$b_0, b_1$ and $b_2$. Then the system can be solved uniquely, producing
a 3-parameter family of polynomials with 4 rational factors
which we do not write down here. We have been unable to find additional
examples in this family for which we could verify Beilinson's conjecture,
either because the coefficients were too large or because
we did not obtain enough elements in $\INT$ of the corresponding curve $C$.
Besides, there is a somewhat neater construction, which gives
more rational roots:

\ExampleNamed{$\mathbf{d=5}$}
\label{ex5_5}
In degree 5 it is in fact possible to describe the polynomials of the
desired form that factor completely over the rationals as follows.
Recall from Example \ref{ex5_2}
that every rational root of $t(x)$ is of the form $-r^2$ for some
rational~$r$, so if $ t(x) $ factors completely then we must have
$$
  t(x) = x^5 + (b_2 x^2 + b_1 x + b_0)^2 = \prod_{i=1}^5(x+r_i^2) \>.
$$
Assume that the $|r_i|$ are pairwise distinct and non-zero
(otherwise $ t(x) $ has a double root).  The two formulas for $t(x)$ give
two factorizations of $-t(-x^2)$, and hence
$$
  (x^5 + b_2 x^4 - b_1 x^2 + b_0)(x^5 - b_2 x^4 + b_1 x^2 - b_0)
    \>\>=\>\> \prod_{i=1}^5(x-r_i)\prod_{i=1}^5(x+r_i) \>.
$$
Call the two quintics on the left $q(x)$ and $-q(-x)$ respectively.
Clearly in every pair $\pm r_i$ one of the numbers is a root of
$q(x)$ and the other a root of $q(-x)$. By changing the signs
of some of the $r_i$ if necessary, we can assume that all of
them are roots of $q(x)$. As the quintic $q(x)$ has no $x^3$ and
$x$ terms, Newton's formulae imply that the numbers $s_i\!=\!r_i^{-1}$
satisfy the two equations
\beq\label{cubsur}
  \sum_{i=1}^5 s_i = 0 \quad  \text {and} \quad \sum_{i=1}^5 s_i^3 = 0 \>.
\eeq
Conversely, any five-tuple $\{s_i\}$ of non-zero rationals satisfying (\ref{cubsur})
and such that $s_i\ne\pm s_j$ for $i\ne j$
gives rise to $t(x)=\prod(x+s_i^{-2})$ of the desired form.
Note that the tuple $\{-s_i\}$ gives the same polynomial.

Now in order to find rational solutions of (\ref{cubsur}), note that
these equations define a cubic surface in $\P^4$. So given two
rational solutions to (\ref{cubsur}), the line joining them intersects the
surface in a third point which is also rational.
Thus, starting with some obvious solutions such as $(a,-a,b,-b,0)$ and
its permutations or some other experimentally found small solutions like
\beq\label{cubsurex}
  (1, 5,-7,-8, 9), \quad (2,4,-7,-9,10), \quad (1,14,-17,-18, 20),
  \quad\ldots
\eeq
one can construct as many other solutions as one wants by using this
``chord construction'' and by scaling and permuting the coordinates.

\begin{remark}
The surface defined by~\eqref{cubsur} is the famous cubic surface
of Clebsch and Klein, which also occurs as a Hilbert modular surface and was studied in detail in 
\cite{Hir}.  It might be worth investigating why this surface shows up in this context.
\end{remark}

\ExampleNamed{$\mathbf{d}$ even}
\label{exd}
Finally, here is a method to produce examples for even $d=2m$. First, if
\beq\label{c2eq1}
  t(x) = -c^2x^{2m} + f(x)^2/4 = (cx^m + f(x)/2)(f(x)/2-cx^m)
\eeq
with $ f(x) = 2 c x^m + b_gx^g + \dots +b_0 $ and $ cb_g\neq 0 $,
then we could take the first factor to be of the form $2c(x-a_i)\dots(x-a_m) $,
which gives us at least $ m+1 $ rational factors of $ t(x) $.
In order to get more rational factors, it is easier to put
$ \tilde f(x) = x^m f(1/x) $, so that, with $g(x)=\tilde f(x)/2+c$ and $k=2c$, 
\beq\label{c2eq2}
  x^{2m}t(1/x) = \tilde f(x)^2/4 - c^2 
= g(x)(g(x)-k) \>,
\eeq
and the problem simply becomes to find polynomials $g(x)$
and non-zero constants $k$ such that $g(x)$ and $g(x)-k$
have no multiple roots and a lot of rational factors.
Strictly speaking we should also ensure that
(\ref{c2eq2}) has no constant term and a non-zero linear term, so that (\ref{c2eq1})
is of degree $2m-1$, a condition which has to be satisfied to
get the correct shape of $t(x)$.
But this 
can be always achieved by a translation if there is a rational
root.

Scaling $ g(x) $ and $ k $ does not change our problem, so we
can assume $ g(x) $ has leading coefficient 1.
If we choose $g(x) = (x-b_1)\dots(x-b_m)$ we get a family of good polynomials with at
least $m+1$ rational factors, having $b_1,\dots,b_m$ and $ k $ as parameters.
Now the question is how to improve this by forcing a polynomial of the form
\beq*
  g(x) - k = (x-b_1)\cdots(x-b_m) - k
\eeq*%
to have more rational factors.
If we choose $k=\prod(-b_i)$, then $g(x)-k$ has a factor $x$,
so we get a family with at least $m+2$ rational factors, namely
$x, x-b_1,\dots,x-b_m$ and the remaining factor of degree $m-1$.

For $d=6$ $(m=3)$
this remaining factor is quadratic.
So with a suitable
rational parametrization one can construct a ``universal'' family of
polynomials which factor completely over the rationals.
One example of such a family is given by
$$
\begin{array}{rcl}
  g(x) & = & (x-1)(x-r-s)(x-rs), \cr
  k & =  & rs(r-1)(s-1), \cr
  g(x)-k & = & (x-r)(x-s)(x-rs-1) \>.
\end{array}
$$

For larger degrees, examples can be found as follows.
Without loss of generality we may assume that $b_1,\dots,b_m$
are integers. Choose bounds $B, C$ and search through all
integers $-B\le b_1<\dots<b_m\le B$.  Let $g(x)=(x-b_1)\dots(x-b_m)$
and compute $g(-C),\dots,g(C-1),g(C)$. If a number $k\ne 0$ occurs
more than once in this list of values, then $g(x)-k$ has several integral roots
in the range from $-C$ to $C$, each of which yields one linear factor.
This gives a simple method to look for polynomials~\eqref{c2eq2}
with even more rational roots.


Here are a few examples with $m\!=\!4$ and $m\!=\!5$
where both $g(x)$ and $g(x)-k$ factor completely, sorted according to $k$.
The columns in the two tables below contain $k$, the roots of $g(x)$ and
the roots of $g(x)-k$ respectively.
\def\entry#1|#2#3{\footnotesize$#1$ & \footnotesize$#2$ & \footnotesize$#3$\cr}
\def\enthead#1{
  \begingroup
  \def\arraystretch{0.85}
  \begin{tabular}{|l|l|l|}
  \hline
  \multicolumn{3}{|c|}{\footnotesize $m=#1$}\cr
  \hline
  \footnotesize$\hfill k\hfill$ & \footnotesize$\hfill g(x)=0 \hfill$ & \footnotesize$\hfill g(x)=k \hfill$\cr
  \hline
}
\def\enttail{\hline\end{tabular}\endgroup}
$$
\enthead4
\entry      2520|{9,7,-4,-10}{11,2,0,-11}
\entry      7560|{13,12,-5,-14}{16,7,-2,-15}
\entry     10080|{14,12,-4,-15}{17,6,0,-16}
\entry     10080|{17,14,-9,-17}{19,11,-7,-18}
\entry     12600|{17,12,-6,-17}{19,8,-3,-18}
\entry     13860|{15,12,-4,-17}{18,5,1,-18}
\entry     15840|{15,14,-3,-17}{19,5,3,-18}
\entry     18720|{19,12,-5,-19}{21,7,-1,-20}
\entry     25200|{19,14,-3,-20}{22,5,4,-21}
\entry     27720|{21,12,-5,-22}{23,6,0,-23}
\enttail
\quad
\enthead5
\entry     50400|{12,11,-1,-6,-14}{14,6,4,-9,-13}
\entry     50400|{14,10,1,-9,-15}{15,6,5,-11,-14}
\entry    110880|{15,13,-3,-4,-16}{17,8,4,-9,-15}
\entry    110880|{16,10,-3,-7,-16}{17,5,4,-12,-14}
\entry    272160|{17,16,-1,-7,-20}{20,8,7,-11,-19}
\entry    327600|{20,11,-5,-9,-21}{21,5,4,-15,-19}
\entry    393120|{23,17,-6,-15,-24}{24,15,-3,-19,-22}
\entry    554400|{21,19,2,-11,-24}{24,11,9,-14,-23}
\entry    554400|{22,17,-1,-9,-24}{24,11,6,-13,-23}
\entry   1058400|{23,22,-1,-8,-27}{27,13,8,-13,-26}
\enttail
$$
\smallskip

One should note here that such a search produces many more
examples than just the ones above. The tables actually start
$$
\enthead4
\entry       180|{5,4,-3,-4}{6,2,-1,-5}
\entry       360|{6,4,-3,-5}{7,1,0,-6}
\entry       504|{7,5,-4,-6}{8,3,-2,-7}
\entry       720|{7,4,-4,-7}{8,1,-1,-8}
\entry    \ldots|\ldots\ldots
\enttail
\qquad
\enthead5
\entry      5040|{8,7,-1,-5,-9}{9,5,1,-7,-8}
\entry     10080|{9,7,-2,-4,-10}{10,4,2,-7,-9}
\entry     50400|{12,11,-1,-6,-14}{14,6,4,-9,-13}
\entry     50400|{14,10,1,-9,-15}{15,6,5,-11,-14}
\entry    \ldots|\ldots\ldots
\enttail
$$
However, most of the examples are induced in the sense that
modulo a linear substitution, $p(x)=g(x)(g(x)-k)$ is of the form
$q(x^2)$ for some $q(x)$ of degree $m$.
For instance, for the first entry for $m\!=\!4$ ($k\!=\!180$) we have
$$
  p(x) = (x-5)(x-4)(x+3)(x+4)\cdot(x-6)(x-2)(x+1)(x+5)
$$
and, after rearranging the factors, we find that
$$
\begin{array}{rl}
 p(-x\!+\!\frac12)
 =&\!\!\!\!  (x\!-\!\frac{11}2)(x\!-\!\frac92)(x\!-\!\frac72)(x\!-\!\frac32)(x\!+\!\frac32)(x\!+\!\frac72)(x\!+\!\frac92)(x\!+\!\frac{11}2)\cr
 =&\!\!\!\!  x^8 - 65x^6 + \frac{10979}{8} x^4 - \frac{164385}{16} x^2 + \frac{4322241}{256} \>.
  \phantom{x^{5 \choose 5}} 
 \end{array}
$$
The curve $y^2=t(x)$ is therefore a cover of an elliptic curve
$y^2=x^4-65x^3+\cdots$. In the same way, the first two lines
for $m\!=\!5$ simply come from genus 2 curves, which we have already
seen in the cubic surface construction: they are given
by the first two 5-tuples in (\ref{cubsurex}).


Finally, for $m=6$ $(g\!=\!5, d\!=\!12)$ one can easily find several non-trivial
examples where $g(x)(g(x)-k)$ factors completely except for one quadratic
factor. It is harder to find non-trivial examples which factor completely,
but they do exist. In the notation of (\ref{c2eq2}) the smallest one
is
$$
\begin{array}{rl}
 &\!\!\!\! (x^6\!+\!2x^5\!-\!787x^4\!-\!188x^3\!+\!150012x^2\!-\!149040x\!-\!3326400)^2 \>-\> 3326400^2\cr
 =&\!\!\!\! (x\!-\!22)(x\!-\!20)(x\!-\!18)(x\!-\!12)(x\!-\!10)(x\!-\!1)x(x\!+\!7)(x\!+\!15)(x\!+\!18)(x\!+\!23)(x\!+\!24)\>.
\end{array}
$$


\section{Integrality of the elements} 

Construction \ref{maincons} shows how to construct
hyperelliptic curves $C/\Q$ with elements in $\KC/\torsion$.
However, we are not yet in a position to compute regulators and
test Beilinson's conjecture~\ref{beilinsonconjecture},
because we still need our elements to be {\em integral\/},
that is, to lie in $\INT\subseteq \KC/\torsion$.  (See~\eqref{INTdef}.)

Let $\CZ/\Z$ be a regular proper model of $C/\Q$. 
Recall that an element $\alpha$ of $K_2(\Q(C))$
lies in $K_2^T(C)$ if the tame symbol $T\iP(\alpha)$ is trivial
for all $P\in C(\Qbar)$. Recall also that $\alpha$ lies in
$ K_2^T(\calC) $ if the tame symbol is trivial for each irreducible
curve $\D$ in $\CZ$; see~\eqref{KTcalC} or~\eqref{realtame}. This means that
apart from being trivial for 
\lq\lq horizontal\rq\rq\ curves (which come from points $P\in C(\Qbar)$),
the tame symbol of $\alpha$ must be trivial for all irreducible components of
the fibers $ \CZ_p $ of $\CZ\to\Spec\Z$.
This gives additional conditions on $\alpha$ for all primes $p$ of $\Z$.
(One can show though that, to a given  $ \a $, one can always
add an element in $ K_2(\Q) $ such that the sum
satisfies this condition for each prime $ p $ for which the fiber $\calC_p$ is smooth over $ \F_p $.)


So, in general, given a curve $C/\Q$ and $\alpha\in\KC$ the way to verify
that $\alpha$ gives rise to an element of $ \INT$ is to find a regular model $\CZ/\Z$
and then check that the tame symbol of $ \alpha $ at each ``vertical'' curve on $\CZ$
is trivial.  In practice, finding such a model means
starting with any equation of $C$ with integer coefficients,
which defines an arithmetic surface,
and performing blow-ups until we obtain a regular surface.
This, however, is a complicated
process, so we will try to deduce integrality from the behaviour
of $\alpha$ on the original (possibly singular) arithmetic surface.

This can be done in fair generality for the families of examples
that we are interested in.  We repeat our notation for the sake
of easy reference, so we consider a curve of genus
$ g $ as in Construction \ref{maincons},
defined by
\beq\label{eqdeffrepeat}
  y^2 + f(x) y + x^d = 0
\eeq
in one of the two cases
\beq\label{eqdeff2repeat}
\begin{array}{llll}
\qquad\text{(a)} &  d=2g+1,  & f(x)= b_g x^g + \ldots + b_1 x + b_0 \>; \cr
\qquad\text{(b)} &  d=2g+2,  & f(x)= 2x^{g+1}+b_g x^g + \ldots + b_1 x + b_0\text{ and } b_g \neq 0\>, \cr
\end{array}
\eeq
where $ t(x) = -x^d + f(x)^2/4 $ has no multiple roots, so that
in particular $ b_0 \neq 0 $.
Assume further that $f(x)$ has integer coefficients, so~\eqref{eqdeffrepeat}
defines an arithmetic surface over~$\Z$.

\begin{theorem}\label{thmintegrality}
Let $C/\Q$ be defined by~\eqref{eqdeffrepeat} and~\eqref{eqdeff2repeat},
where $b_0,\dots,b_g$ are integers.
Let $m(x)\in\Z[x]$ be a non-constant factor, irreducible over $ \Z $,
of the 2-torsion polynomial $t(x)\!=\! -x^d\!+\!f(x)^2/4$
and let $M$ be the class of
$$
  \Bigl\{\>\frac{y^2}{x^d}\>,\>\frac{m(x)}{m(0)}\>\Bigr\} \>,
$$
in \KCtors.
We have
\begin{enumerate}
  \item\label{int}    if $m(0)=\pm 1$ then $M\in\INT$;
  \item\label{nonint} if there is a prime dividing  $m(0)$ but not every $b_i $ then $nM\not\in\INT$ for any $n\ne 0$.
\end{enumerate}
Moreover,
\begin{enumerate}
  \item[(3)]
        if $ b_0 = \pm 1 $  and $ \specelt $ is as in Proposition~\ref{Vproposition} then
        $ \specelt $ is in $ \INT $ if $ d=2g+1 $ and $ 2 \specelt $ is in $ \INT $ if $ d=2g+2 $.
\end{enumerate}
\end{theorem}

\begin{proof}
Denote
$$
  \tM=\{f_1,f_2\}, \qquad
  f_1 = \frac{y^2}{x^d}, \quad f_2 = \frac{m(x)}{m(0)} \>;
$$
then $ \tM $ is an element of $ K_2^T(C) $ (see~\eqref{normcalculation}) and $ M $ is its class in $ K_2^T(C)/\torsion $.
We shall show that, for a specific regular proper model $ \calC $ of $ C $, 
$ \tM\in K_2^T(\calC) $ if $ m(0)=\pm1 $, so that $ M \in \INT $.
But, for the same model, we shall show that 
under the asssumptions of~\eqref{nonint}, 
$ T_\D (\tM) $ is not torsion  for some irreducible component $ \D $ of the fiber $ \calC_p $.
Therefore no non-zero multiple of $ M $ can
lie in $ \INT \subseteq K_2^T(C)/\torsion $, because, up to torsion,
$ K_2^T(\calC) $ is a subgroup of $ K_2^T(C) $ that does not depend on the choice of $ \calC $.
Finally, for~(3) we show that $ \{-b_0y,-x\} $ or $ 2\{-b_0y,-x\} $ is in $ K_2^T(\CZ) $.

Before we can give the proof proper, we need some preliminaries.

Let $p$ be a prime and $F$ the fiber above $p$ of the arithmetic surface defined
by~\eqref{eqdeffrepeat}.
It is clear from the equation that any irreducible component $D$ of $F$ dominates the $x$-axis and that $ x $, $ y $ and
$ m(x) $ define rational functions on $ D $ that do not vanish identically on it.

We want to know that if $ p $ does not divide every $ b_i $
then the function $y^2/x^d$ is non-constant on every such $D$.
If it were constant then $y^2=kx^d$ on $D$ for some
constant $k\neq 0$.  So $f(x)y = -(k+1)x^d$ on $ D $
and squaring this equation yields $ f(x)^2 kx^d = (k+1)^2x^{2d} $.
Since $D$ dominates the $x$-axis, this is only possible if this is an
identity of polynomials, which implies that $b_g x^g +\dots+ b_0=0\mod p $,
contradicting the assumption on the $b_i$.

Let $\CZ/\Z$ be the regular model of $C/\Q$ obtained from
the arithmetic surface defined by~\eqref{eqdeffrepeat} by repeatedly blowing up the
singularities.
If $\D\subset\CZ$ is an irreducible curve mapping onto $D$, then
the function $f_1 = \frac{y^2}{x^d}$ is non-constant along $\D$ as well.

We are now ready to prove~\eqref{nonint}.
Let $p$ divide $ m(0) $ but not every $ b_i $.
Let $D$ and $\D$ be as above, so that $v\si\D(f_1)=0$ and $v\si\D(f_2)<0$.  Then the tame symbol $T\si\D(\tM)$
is a non-zero power of $f_1$ and is therefore non-constant
and hence non-torsion.

We now move to the proof of~\eqref{int} so we assume that $m(0)=\pm 1$.
Then for any prime $p$ and irreducible component $\D$ of $ \CZ_p $ surjecting onto a component $ D $
of the fiber above $ p $ of~\eqref{eqdeffrepeat},
$v\si\D(f_1)=v\si\D(f_2)=0$, hence $T\si\D(\tM)\!=\!1$.

However, this does not prove that $v\si\D(f_1)=v\si\D(f_2)=0$
for {\em every\/} component $\D$ of the fiber of $\CZ\to\Spec\Z$ above $p$:
such a $ \D $ could also map to a singular point of the model~\eqref{eqdeffrepeat}.
In this case if, say, $f_1$ happens to have a zero or a pole passing just through
this point, then it can happen that $v\si\D(f_1)\ne0$.

We claim that at every singular point $P$ of the surface defined
by~\eqref{eqdeffrepeat} either $f_1$ or $f_2$ is regular and equal to~1. 
This implies that $T_\D(\tM) = 1 $ for every irreducible curve $\D$ of $\CZ$ mapping to $P$, proving~\eqref{int}.

In fact, our claim holds for all singular points in the fiber $ F $ of~\eqref{eqdeffrepeat} above
a prime $ p $, among which are those singularities of $ \CZ $ that lie above $ p $.
In order to see this, first consider the point at infinity in $ F $.
The arithmetic surface defined by~\eqref{eqdeffrepeat} has a chart at
infinity which can be obtained by letting $x=1/\tilde x$ and
$y=\tilde y/{\tilde x}^{g+1}$. So for $d=2g+1$, the
equation at infinity is
$$
  \tilde y^2 + (b_g \tilde x + \ldots + b_1 \tilde x^g + b_0 \tilde x^{g+1})
  \,\tilde y + \tilde x = 0
$$
and the point at infinity in $ F $ ($\tilde x = \tilde y=0$ mod $p$) is non-singular.
If $d=2g+2$, then the equation is
$$
  \tilde y^2 + (2 + b_g \tilde x + \ldots + b_1 \tilde x^g + b_0 \tilde x^{g+1})
  \,\tilde y + 1 = 0
$$
and the function $f_1=\tilde y^2$ is regular and equal to 1 at the point at infinity,
($\tilde x, \tilde y) = (0,-1)$ mod~$p$.

For finite $x$, if $P=(x_0,y_0)\mod p$ is a singularity of $ F $
then $f(x_0)=-2y_0$. Then either $x_0\ne 0$, in which case $ f_1 $
is regular and equal to 1 at $P$ from the equation of the curve,
or $x_0=0$, in which case $f_2$ is regular and equal to 1 at
$ P $ because $p$ does not divide $m(0)$. This completes the proof
of~\eqref{int}.

For~(3), assume that $ b_0 = \pm 1 $, so that $ \specelt $ is the class in $ \KC/\torsion $
of the element $ \{-b_0y,-x\} $ in $ \KC $.
As above, there can only be a problem at a singular point of the arithmetic surface defined by~\eqref{eqdeffrepeat},
and that point is a singularity on the fiber $ F $ above some prime $ p $.
But if $ P = (x_0, y_0) $ modulo~$ p $ is any finite singularity of $ F $,
then $ 2y_0+f(x_0)=0 $, hence $ y_0^2=x_0^d $.  We cannot have $ (x_0,y_0) = (0,0) $ since this would give $ \pm 1 = b_0 = f(x_0) = -2 y_0 = 0 $.
Therefore neither $ -b_0y = \mp y $ nor $ -x $ can have a zero or a pole passing through $ P $, showing that
those functions have trivial valuation along any irreducible curve $ \D $ in  $ \CZ $ mapping to $ P $.
For $ d = 2g+1 $ we saw before that the point at infinity of $ F $ is non-singular,
so that $ \specelt $ is in $ \INT $.  For $ d=2g+2 $ we consider
instead the class $ 2\specelt $ of $ 2 \{-b_0y,-x\} = \{ y^2, x \} = \{ y^2/x^{2g+2}, x \} = \{ \tilde y^2, 1/\tilde x \} $,
and at the point at infinity of $ F $,  $ (\tilde x, \tilde y) = (0, -1) $
modulo~$ p $, $ \tilde y^2 $ is regular and equal to~1.
\end{proof}

\begin{remark}
We cannot drop the assumption in~\eqref{nonint} of Theorem~\ref{thmintegrality} that the prime does not divide each $ b_i $,
i.e., the unrestricted converse of~\eqref{int} does not hold.  
In fact, given $ m(x) $ as in the theorem with $ m(0) = \pm1 $, so that the corresponding $M $ is in $ \INT$, fix a prime
$ p $ and a positive integer $ s $.  Then the substitution $ (x,y) \mapsto (p^{-2s}x,p^{-ds}y)$.
transforms $C$ into an isomorphic curve $ \tC $ defined by
\[
  y^2 + \tilde f(x) y + x^d = 0
\>,
\]
where $ \tilde f(x) = p^{ds} f(p^{-2s}x) $ so that all $ \tilde b_i $ are divisible by $ p $.
Choosing $ k \geq 0 $ minimally such that $ \widetilde m(x) = p^k m(p^{-2s}x) $ is in $ \Z[x] $,
$ \widetilde m(x) $ is irreducible in $ \Z[x] $ and is a non-constant
factor of the 2-torsion polynomial $  -x^d + \tilde f(x)^2/4 = p^{2ds} t(p^{-2s}x) $.
We can ensure that $ k>0 $ by choosing $ s $ sufficiently large, so that $ p $ divides $ \widetilde m(0) $.
Then the class of $ \{y^2/x^d,\widetilde m(x)/ \widetilde m(0) \} $ is in $ K_2(\tC;\Z) $ since it corresponds
to $ M $, but $ p $ divides $\widetilde m(0) $ as well as each $ \tilde b_i $.
\end{remark}


\section{Computing the Beilinson regulator} 
\label{regulat}

Let $ C $ be a hyperelliptic curve of genus $ g $ over $ \C $,
so the map $ \phi : C \to \P_\C^1 $ has $ n = 2 g + 2 $ points of ramification,
$ P_1,\dots,P_n $, and let $ X = C(\C) $ be the associated Riemann surface.
It is not difficult to check that a basis of $
\hsin 1 (X; \Z) $ consists of liftings to $ X $ of simple loops in $
\P_\C^1 $ around exactly $ P_i $ and $ P_{i+1} $, for $ i = 1,\dots,
2g $.  If a model of $ C $ is defined by an equation $ y^2 = f(x) $
with $ f(x) $ in $ \C[x] $ of degree $ 2 g + 1 $ without multiple
roots, the point in $ C $ above the point at infinity in this model
will be among the ramification points, and the other ones will be the
points above the roots of $ f(x) $ in $ \C \subset \P_\C^1 $.  We therefore work with
simple loops in $ \C $ around exactly two roots of $ f(x) $, not passing through any root of $ f(x) $, and lift them to loops
on~$ X $.

\begin{figure}[htb]
\epsfig{file=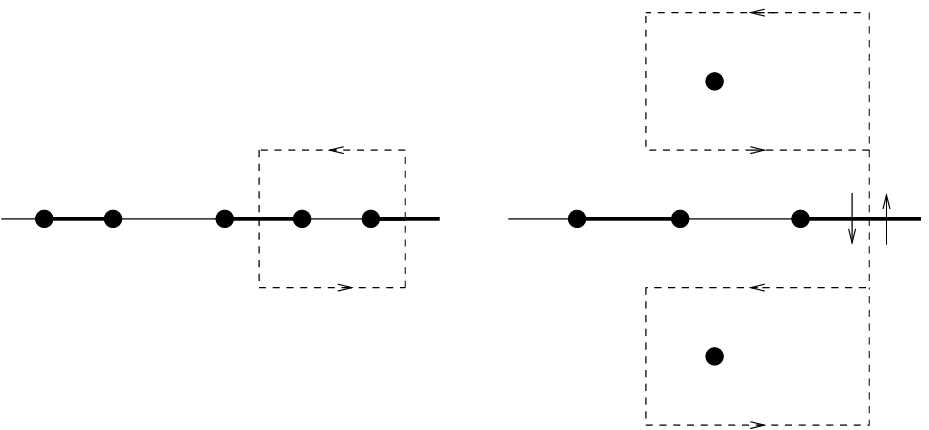}
\caption{}
\label{paths}
\end{figure}

In our case $ f(x) $ belongs to $ \Q[x] $ and hence to $ \R[x] $,
and we want to keep track of the action of complex conjugation
on $ X $.  We will discuss the case when the leading coefficient
of $ f(x) $ is positive, which we can always achieve by replacing
$ x $ with $ -x $ if necessary.
The set $ C(\R) $ of real points on $ X $ consists of $ \infty_C $ together
with the points above those $ x $ in $ \R \subset \C $ 
where $ f(x) \geq 0 $.   In Figure~\ref{paths}, the latter is indicated
by the thick part of $ \R $ in $ \C $, the thick
dots being the
roots of $ f(x) $.   

If $ f(x) $ has only real roots $ P_1 < P_2 < \dots < P_{2g+1} $,
then let  $ \c_i $ be
a lift of a loop around $ P_{2i} $ and $ P_{2i+1} $ ($ i = 1,\dots,g $)
(this is illustrated on the left in Figure~\ref{paths}),
and let $ \delta_i $ be a lift of a loop around $ P_{2i-1} $ and $ P_{2i} $
($ i = 1,\dots,g $).
Since all of the $ \c_i $ have a fixed point
under the action of complex conjugation, which lies
above the intersection with the thick part of the real line,
it is easy to check that the $ \c_i $ lie in $ \hsin 1 (X; \Z)^- $
because complex conjugation reverses the orientation.
Also, by shrinking the $ \delta_i $ as much as possible, so that they lift to a subset of $ C(\R) $,
it is easy to see that  the $ \delta_i $ lie in $ \hsin 1 (X; \Z)^+ $.
So in this case $ \hsin 1 (X; \Z) $ decomposes as a direct sum
$ \hsin 1 (X; \Z)^- \oplus \hsin 1 (X; \Z)^+ $, the $ \c_i $
form a basis of $ \hsin 1 (X; \Z)^- $, and the $ \delta_i $ form a basis of
$ \hsin 1 (X; \Z)^+ $.

If $ f(x) $ does not have only real roots, let $ P_1,\dots, P_{2m+1} $
be the real roots of $ f(x) $, and let $ Q_1, \overline{Q_1} ,\dots, Q_{g-m}, \overline{Q_{g-m}} $
be its  non-real roots, as pairs of complex conjugated numbers
with $ \text{Im}(Q_j) > 0 $.
We now define $ \c_1,\dots,\c_m $ as lifts of simple loops around $ P_2 $ and
$ P_3 $, $ P_4 $ and $ P_5 $, $ \dots $, $ P_{2m} $ and $ P_{2m+1} $,
and $ \c_{m+1},\dots,\c_g $ as lifts of simple loops around $ Q_1 $ and $ \overline{Q_1} $, $ Q_2 $
and $ \overline{Q_2} $, etc., intersecting the thick part of the
real axis as illustrated on the right in Figure~\ref{paths}.
Then it is easy to check that $ \c_1,\dots,\c_g $ are in $ \hsin 1 (X; \Z)^- $.
Also, if $ \delta_1,\dots,\delta_g $ are lifts to $ X $ of simple loops around $ P_1 $ and $ P_2 $,
$ \dots $, $ P_{2m-1} $ and $ P_{2m} $, as well as around
$ P_{2m+1} $ and $ Q_1 $,
$ Q_1 $ and $ Q_2 $, $ Q_2 $ and $ Q_3 $, etc., then the $ \delta_j $
complement the $ \gamma_j $ to a basis of $ \hsin 1 (X; \Z) $.
Note that in Figure~\ref{conj}
\begin{figure}[htbp]
\begin{center}
\begin{picture}(0,0)%
\includegraphics{conj.pstex}%
\end{picture}%
\setlength{\unitlength}{2171sp}%
\begingroup\makeatletter\ifx\SetFigFont\undefined%
\gdef\SetFigFont#1#2#3#4#5{%
  \reset@font\fontsize{#1}{#2pt}%
  \fontfamily{#3}\fontseries{#4}\fontshape{#5}%
  \selectfont}%
\fi\endgroup%
\begin{picture}(6645,3684)(5989,-4603)
\put(8476,-1561){\makebox(0,0)[lb]{\smash{{\SetFigFont{10}{12.0}{\rmdefault}{\mddefault}{\updefault}{$\delta$}%
}}}}
\put(9301,-1486){\makebox(0,0)[lb]{\smash{{\SetFigFont{10}{12.0}{\rmdefault}{\mddefault}{\updefault}{$\gamma$}%
}}}}
\put(6226,-3361){\makebox(0,0)[lb]{\smash{{\SetFigFont{10}{12.0}{\rmdefault}{\mddefault}{\updefault}{$\delta'$}%
}}}}
\end{picture}%

\caption{}
\label{conj}
\end{center}
\end{figure}
we still have $ \delta' \circ \s = \delta' $
as in the case when $ f(x) $ has only real roots, but
we do not get a splitting of $ \hsin 1 (X; \Z) $ into $ + $ and
$ - $-parts.  Namely, if $ \delta $ is the lift of a loop around $ P_{2m+1} $
and $ Q_1 $, and $ \c $ is the lift of a loop around $ Q_1 $
and $ \overline{Q_1} $ passing through the thick part of the
real line to the right of $ P_{2m+1} $, then $ \delta - \delta \circ \s = \c $, with $ \c $
in $ \hsin 1 (X; \Z)^- $, provided $ \delta $ and $ \c $ are
chosen in such a way that they are in the same branch of $ C(\R) $
above the thick part of the real line through which both loops
in $ \C $ pass.
(An illustration of this is given on the right in the left part
of Figure~\ref{conj}.)
This holds because looping twice around a single point lifts to the trivial element in $ \hsin 1 (X; \Z) $.
If we could write $ \delta = \delta^+ + \delta^- $ with $ \delta^\pm $
in $ \hsin 1 (X; \Z)^\pm $, then this would lead to $ 2 \delta^- = \delta - \delta\circ\s = \c $,
which is not possible since $ \c $ is part of a basis of $ \hsin 1 (X;\Z) $.
In summary we find in this case that $ \{\c_1,\dots,\c_g \} $ is
a basis of $ \hsin 1 (X; \Z)^- $, and that $ \{\delta_1,\dots,\delta_g \} $
is a basis of $ \hsin 1 (X; \Z)/\hsin 1 (X; \Z)^- $.

For practical purposes, we choose all our loops in $ \C $ as concatenations
of line segments, which makes computing and parametrizing their lifts to $ X $
using the analytic continuation of square roots particularly easy.
Apart from the two types already described (lifts of a simple loop
around two consecutive real roots of $ f(x) $, or of a simple
loop around two conjugate non-real roots of $ f(x) $ intersecting
the thick part of the real line, both illustrated in Figure~\ref{paths}),
we also use a third type in order to avoid getting close to roots
of $ f(x) $ unnecessarily, thus speeding up the calculations
in numerical integration.
If, for example, $ Q_1 = x_1 + i y_1 $ and $ Q_2 = x_2 + i y_2 $ are two non-real roots of $ f(x) $
which are close together, and $ 0 < y_1 \leq y_2 $, then,
keeping a lift of a loop around the pair $ Q_1, \overline{Q_1} $
meeting the thick part of the real line,
we can replace a lift of a loop around the pair $ Q_2, \overline{Q_2} $
with a loop around both pairs $ Q_1, \overline{Q_1}, Q_2, \overline{Q_2} $,
still meeting the thick part of the real line, as indicated on
the right in Figure~\ref{conj}.  This is just
the sum of (compatible) lifts of a loop around $ Q_1, \overline{Q_1} $
and a loop around $ Q_2, \overline{Q_2} $, both meeting the thick
part of the real line in a common point.
A similar method is used if more than two non-real roots of $ f(x) $
are close together.

Finally, we choose our loops also such that they do not
pass through the image of the other torsion point under the map
$ \phi : X \to \P_\C^1 $. 
This way, we can immediately compute 
$ \langle \c , \{f,g\} \rangle = \frac{1}{2 \pi} \int_\c \eta(f,g)
$, with $ \eta(f,g) $ as in \eqref{eta},
and integrate numerically over the lifts of the paths using the
analytic continuation of $ y = \sqrt{f(x)} $ along the paths.
The method works quite well in practice and we can easily obtain
twelve or fifteen decimals of precision for the regulator.


\section{Examples} 
\label{examples}

In this section we test Beilinson's conjecture \ref{beilinsonconjecture}
as explained in Remark~\ref{untractable}.  We do this for
hyperelliptic curves over $\Q$, of genus 2, 3, 4 and 5.
Using the results of Section \ref{sFam}, we construct curves
of the form \eqref{eqdeffrepeat} and~\eqref{eqdeff2repeat} whose 2-torsion polynomials
have many rational factors. As Construction~\ref{maincons} shows, every
such rational factor gives an element of $\KC/\torsion$.
Theorem~\ref{thmintegrality} shows which of these are in $ \INT $.

Our main sources of examples are the constructions described in 
Examples~\ref{ex5_2} and~\ref{exd} (or variations of them).
Unfortunately, the more sophisticated constructions of
Examples~\ref{ex5_3} and~\ref{ex5_5} tend to produce
either not enough symbols satisfying the integrality condition or
else curves of very high conductor for which  we cannot compute the $L$-value.

\ExampleNamed{genus 2}\label{exgen2}
Taking $m\!=\!1$ and $n\!=\!-2$ in
Example \ref{ex5_2} yields a one-parameter family of curves
$$
  \def\sp{\!+\!}
  y^2 = x^5 +
    \Bigl(
      (k + 3) x^2 + (5 k + 4) x + 4 k
    \Bigr)^2 \>.
$$
Further, if $k$ is an integer which is divisible by 4,
then the curve can be transformed to one of the form
in Theorem \ref{thmintegrality}.
Thus write $k\!=\!-4b$ with $b\in\Z$. The curve is then isomorphic to
$$
  C_b:\quad y^2 + ((4b-3)x^2-(5b-1)x+b) y + x^5 = 0 \>.
$$
For non-zero $b$ this is a non-singular genus 2 curve of discriminant
$$
  \Delta(C_b) = 3^2 b^7 (8b+3)^2 (9b^3-23b^2+4b-1) \>.
$$
Its 2-torsion polynomial has 3 rational factors, namely (up to a constant)
$$
\begin{array}{rcl}
  m_1&=&x-1, \cr
  m_2&=&4x-1, \cr
  m_3&=&x^3-(4b^2\!-\!6b\!+\!1)x^2+(5b^2\!-\!b)x-b^2 \>.
\end{array}
$$
Recall from Construction~\ref{maincons} that
each $m_i$ gives an element $M_i\in \KC/\torsion$.
Intersecting the lattice in $ \KC/\torsion $ which they span with $ \INT $ we obtain a sublattice
of integral elements,
$$
  \Lambda\iM = \langle M_1, M_2, M_3 \rangle \>\>\cap\>\> \INT.
$$
Using Theorem~\ref{thmintegrality} we can determine $\Lambda\iM$ exactly
by inspecting the constant terms of the $m_i$.  We find that
$\Lambda\iM = \langle M_1, M_2, M_3 \rangle$ for $b=\pm1$ and
$\Lambda\iM = \langle M_1, M_2 \rangle$ otherwise.

We can now test Beilinson's conjecture for this family.
We constructed at least
two elements ($M_1$ and $M_2$) of $\Lambda\iM$ and we might {\em generally\/}
expect them to be linearly independent. If this is the case,
$\Lambda\iM$ should be a subgroup of finite index in $\INT$, which is supposed
to have rank 2 by the first part of Beilinson's conjecture
\ref{beilinsonconjecture}.  The second part would then imply
that $R(\Lambda\iM)$, the regulator computed using a basis of $\Lambda\iM$, is a non-zero rational
multiple of $L^*(C_b,0)$, the
leading coefficient of $L(C_b,s)$ at $s\!=\!0$. In summary, we expect
\bit
  \item[(a)] $\rk\iZ\Lambda\iM = 2$;
  \item[(b)] $R(\Lambda\iM) = Q L^*(C_b,0)$ for some $ Q \in \Q^* $.
\eit

Note that the prediction $ \rk \Lambda\iM = 2 $ implies that in the case $b=\pm 1$,
where we have $\Lambda\iM=\langle M_1, M_2, M_3\rangle$,
there should be a linear relation among the $M_i$.
Let us illustrate the latter point in the case $b=-1$. We can compute
the images of the $M_i$ under the regulator pairing~\eqref{INTpairing},
as explained in Section~\ref{regulat},
\beq\label{regpair}
  \alpha \>\longmapsto\> (\langle \gamma_1,\alpha \rangle, \langle \gamma_2,\alpha \rangle)
    \in \R^2,
\eeq
where $\gamma_1$ and $\gamma_2$ form a basis of $\hsin 1 (X; \Z)^-$.
Since we expect the image to be a lattice of rank 2, there should be some
integral linear combination of the $M_i$ for
which the image in $\R^2$ vanishes. Using LLL to look for a small \Z-relation
between the images in $\R^2$, we find that
\begin{equation}\label{Mrelation}
   41 M_1 + 56 M_2 - 44 M_3 \>\>{\buildrel?\over\longmapsto}\>\> (0,0) \>.
\end{equation}
To make the relation more transparent we complete the vector
$(41, 56, -44)$ to a unimodular integral matrix and make the corresponding
change of basis:
$$
  \begin{pmatrix} M_1^* \cr M_2^* \cr M_3^* \end{pmatrix} \>\>=\>\>
  \begin{pmatrix} 0 & 0 & 1 \cr 11 & 15 & 0 \cr 41 & 56 & -44 \end{pmatrix}
  \begin{pmatrix} M_1 \cr M_2 \cr M_3 \end{pmatrix} \>.
$$
The numerical values of the regulator pairing \eqref{regpair} on the $M_i$
and the $M_i^*$ are then given by the table
$$
  \begin{array}{|c|r|r|r|r|r|}
  \hline
  \alpha & \hfill M_1 \hfill & \hfill M_2 \hfill & \hfill M_3 = M_1^* \hfill & \hfill M_2^* \hfill & \hfill M_3^* \hfill \cr
  \hline
  \langle \gamma_1,\alpha \rangle & -0.519837 & -6.055409 & -8.191279  & -96.549352 & 0.000000 \cr
  \langle \gamma_2,\alpha \rangle &  3.243279 & -0.869701 &  1.915254  &  22.630553 & 0.000000 \cr
  \hline
  \end{array}
$$
This strongly suggests that $M_3^*$ is in the kernel of \eqref{regpair},
in which case $M_3^*$ should be zero in $\Lambda\iM$ by Beilinson's conjecture.
The regulator $R(\Lambda\iM)$
is now obtained as the absolute value of the determinant of the $2\times 2$
matrix which consists of the two columns for $M_1^*$ and $M_2^*$,
$$
  R(\Lambda\iM) = \left|
    \begin{array}{rr}\!\! -8.191279  & -96.549352    \cr  1.915254  &  22.630553  \end{array}\right| =
    0.456625 \>.
$$
Numerically, $L^*(C_{-1},0)\approx 0.228312$, so
$L^*(C_{-1},0)/R(\Lambda\iM)\approx 0.500000$, which we ``recognize'' as the
rational number 1/2. (Here the agreement is
to at least 18 digits, which was our working precision in this example.)

Recall from Theorem~\ref{thmintegrality} that for $b=\pm1$ the class $\specelt$ of $\{-by,-x\}$
also lies in $ \INT $ and satisfies
$10\specelt = M_1+M_2+M_3 $  by~\eqref{specelt} of Proposition~\ref{relprop}.
For $b=-1$, we can write this as $ 10\specelt = -175M_1^*+15M_2^*-4M_3^* $, so 
if we assume $ M_3^*=0 $, then 
$2\specelt= -35 M_1^* + 3 M_2^*$ is in $\Lambda\iM$ but $\specelt$ is
not.
In other words, the larger lattice $\Lambda=\langle M_1, M_2, M_3, \specelt\,\rangle$
in $\INT$ contains $\Lambda\iM$ as a sublattice of index 2 and
$$
  L^*(C_{-1},0)/R(\Lambda)=1.000000... \>\>{\buildrel?\over=}\>\> 1 \>.
$$

Table~\ref{tab5} summarizes the computations for the curves $C_b$
with $|b|\!\le\!10$. In the third column we describe the lattice $\Lambda$
that we use,
and the last column contains the quotient of
$L^*(C_b,0)$ and the Beilinson regulator $R(\Lambda)$.  In all cases this
quotient
appears to be a rational number of relatively small height.
For $b=\pm1$ we also include the expected relation between the $M_i$.
(We have not included the universal relation $10\specelt=M_1+M_2+M_3$ because
similar relations would clog up later tables.)

Finally, we remark that the Jacobian of each of the curves appearing
in Table~\ref{tab5},  except for $ b=1 $, is not isogenous to the product of two elliptic
curves, even over $ \Qbar $.  This was verified using the procedure
of \cite[Chapter~14, Section~4]{CF}.
In other words, the elements in $ K $-theory and the $ L $-values
do not come from elliptic curves, and the verification of
Beilinson's conjecture does not reduce to the case of genus~1.
The same holds for all the curves in Table~\ref{tab6} (Example~\ref{26example})
and the curves in Table~\ref{tru2} (Example~\ref{extrue})
except the one with $v=\{2,5\}$.

\Notation
In all of the following examples we keep the same notation as
above. For a curve
$$
  C: y^2+f(x)y+x^d=0
$$
as in~\eqref{eqdeffrepeat} and~\eqref{eqdeff2repeat},
denote by $m_i$ the rational factors of the 2-torsion polynomial
$t(x)\!=\!-x^d\!+\!f(x)^2/4$. The associated elements of $\KC/\torsion$
(see \eqref{Mjdef}) are denoted by $M_i$.  Using $ \specelt $ or $ 2 \specelt $ 
to aim for as large a subgroup of $ \INT $ as possible
(see Theorem~\ref{thmintegrality}) we let
$$
  \Lambda' = \left\{
     \begin{array}{ll}
       \langle M_1,\ldots,M_k \rangle, & f(0) \ne \pm 1,\cr
       \langle M_1,\ldots,M_k,\specelt\,\rangle, & f(0)=\pm1, d=2g+1,\cr
       \langle M_1,\ldots,M_k,2\specelt\,\rangle, & f(0)=\pm1, d=2g+2,\cr
     \end{array}
  \right.
$$
and use the sublattice $\Lambda=\Lambda'\cap\INT$ of integral elements.
In all of the examples below we can determine $\Lambda$ by applying
Theorem~\ref{thmintegrality}. Then we verify numerically that $\Lambda$
has rank equal to the genus of $C$ and that the leading coefficient $L^*(C,0)$
is a rational multiple of the regulator $R(\Lambda)$ of $\Lambda$.
To keep the entries for $ \Lambda $ simple, we did not include the universal relation
between $\specelt$ and the $M_i$ of~\eqref{specelt} of Proposition~\ref{relprop} in the tables.

\ExampleNamed{genus 2}
\label{26example}
We can construct another family of genus 2 curves with enough elements in $\INT$
in a way similar to Example~\ref{exd}, this time using 6-torsion points.
The curves are given by
$$
   y^2 + (2x^3-4bx^2-x+b) y + x^6 = 0\,, \qquad b\in\Z_{>0}.
$$
The 2-torsion polynomial has 4 rational factors,
$$
  m_1 = x-b, \quad
  m_2 = 2x-1, \quad
  m_3 = 2x+1, \quad
  m_4 = 4bx^2+x-b,
$$
so we get 4 elements of $\KC/\torsion$, again denoted by $M_1,M_2,M_3,M_4$.
According to~\eqref{urelation} of Proposition~\ref{relprop}
these satisfy the relation $M_1\!+\!M_2\!+\!M_3\!=\!M_4\!$ so
$M_4$ can be dropped.  Furthermore, $M_2$ and $M_3$ are always integral, and
$M_1$ is for $b\!=\!1$.
The numerical results, and the expected relation between the $M_i$
when $b\!=\!1$, are given in Table~\ref{tab6}.


\ExampleNamed{genus 3}
We can look at hyperelliptic genus 3 curves with a rational 7-torsion point and whose 2-torsion polynomial has
two given rational root, cf.~Example~\ref{ex5_2}.  A special case of this is the
two-parameter family given by
$$
  y^2 + ((4b\!-\!3)x^3-(4a\!+\!5b\!-\!1)x^2+(5a\!+\!b)x-a) y + x^7 = 0 \>.
$$
The curves have a rational 7-torsion point $(0,0)$ and two rational
2-torsion points with $x$-coordinates $1$ and $1/4$.
The 2-torsion polynomial has 3 rational factors,
$$
\begin{array}{rcl}
  m_1&\!\!=\!\!&x-1, \cr
  m_2&\!\!=\!\!&4x-1, \cr
  m_3&\!\!=\!\!&x^5 - (4b^2\!-\!6b\!+\!1)x^4 + (5b^2\!+\!8ab\!-\!2b\!-\!6a)x^3\cr
            &&  \qquad\qquad\qquad\qquad - (b^2\!+\!10ab\!+\!4a^2\!-\!2a)x^2 + (5a^2\!+\!2ab)x - a^2 \>.
\end{array}
$$
For $a=\pm 1$ and $ b \in \Z $ we thus get 3 integral symbols and we expect a relation
with $L^*(C,0)$.
Moreover, for $a\!=\!1, b\!=\!3$ we have a further factorization
$$
  m_3 = m_{3,1}m_{3,2} = (x^2\!-\!3x\!+\!1)(x^3\!-\!16x^2\!+\!8x\!-\!1).
$$
so we get 4 elements of $\INT$ and we expect them to be linearly dependent.
The numerical results, and the expected relation for $ a=1 $, $ b=3 $, are summarized in Table~\ref{tab7}.


\ExampleNamed{genus 4, 5}
We can also give some sporadic examples of curves of genus $g\!=\!4$
and $g\!=\!5$, given by
$$
  y^2 + f(x) y + x^{10} = 0
$$
with $f(x)$ of degree 5 and
$$
  y^2 + f(x) y + x^{12} = 0
$$
with $f(x)$ of degree 6 respectively.
These can be found by a brute search for good polynomials with enough
factors in $ \Z[x] $ with constant term $\pm 1$.  In all cases that we found,
we have exactly $g+1$ such factors and the lattice
$\Lambda=\langle M_1,\ldots,M_g,2\specelt\rangle$ is of rank $g$.
The examples and the corresponding numerical results are given
in Tables~\ref{tab10} and~\ref{tab12}.

\ExampleNamed{arbitrary genus}\label{extrue}
Finally, we give examples of curves of arbitrary genus $ g $ with $g$ integral symbols.
Consider a curve of the form
\beq\label{rueq}
  y^2 + \Bigl(2x^{g+1} \pm \prod_{j=1}^g (v_jx+1)\Bigr)\>y + x^{2g+2} = 0 \>,
\eeq
where $v_1<\dots<v_g$ are distinct non-zero integers, and we
are assuming that the 2-torsion polynomial $t(x)$ has no multiple roots.
Then $t(x)$ has at least $g+1$ factors in $ \Z[x] $,
\beq\label{rufac}
  v_1 x+1, \> v_2x+1, \ldots, v_gx+1, \> 4x^{g+1} \pm \prod_{j=1}^g (v_jx+1)
\> .
\eeq
(This is similar to, but different from, the discussion in Example~\ref{exd}
because now,
in $ 4t(x) = -4x^{2g+2} + f(x)^2 = (f(x)-2x^{g+1})(2x^{g+1}+f(x)) $,
we decompose the first factor, which has degree $ g $, rather than the second factor, which has degree
$ g+1 $, completely into linear factors.)
The elements $M_j$ associated to all irreducible factors
in $ \Z[x] $ are integral by Theorem~\ref{thmintegrality},
so we get  at least $g+1$ elements of $\INT$, but we have to take~\eqref{urelation} of Proposition~\ref{relprop} into account.

An equation of the form \eqref{rueq} is not unique for a
given curve.  For even $g$, we may (and will) assume that the sign `$\pm$' in \eqref{rueq}
is `$+$', for we can otherwise replace $ (x,y,\{v_j\}) $ with $ (-x,-y,\{-v_{g+1-j}\}) $.
On the other hand, for odd $g$ the map $(x,y) \mapsto (-x,y) $
gives an isomorphism between the two curves with the same sign associated to
$\{v_j\}$ and $\{-v_{g+1-j}\}$.

Also, if we look at the model at infinity by letting
$y\mapsto y/x^{g+1}$ and $x\mapsto 1/x$, the equation becomes
$$
  y^2 + \Bigl(2 \pm x\,\prod_{j=1}^g (x+v_j)\Bigr)\>y + 1 = 0 \>.
$$
Now it is clear that translating $x$ by $-v_j$ for some $ j $ gives an
equation of the same shape. In other words, $\{v_j\}$ and $\{w_j\}$
yield an isomorphic curve whenever $\{v_j\}\cup\{0\}$ is a translate of
$\{w_j\}\cup\{0\}$. Thus, after translating by $-v_1$
in case $ v_1<0 $, we can assume that all $v_j$ are positive.
Combining this with the above, we see that for odd $ g $,
$ 0 < v_1< \dots<v_g $ and $ 0 < v_g\!-\!v_{g-1} < \dots < v_g\!-\!v_1 < v_g $ give
isomorphic curves, and we have chosen the lexicographically smaller
representative in Tables~\ref{tru3plus} and~\ref{tru3minus} for genus 3 curves.

We can also look at the cases where the last rational factor in \eqref{rufac}
is reducible, for instance, when it has a linear factor $ax-1$ with
$a\in\Z$. Then $1/a$ is a root of this polynomial, so
$$
  \prod_{j=1}^g (v_j/a+1) = \mp 4a^{-g-1}.
$$
It follows that $a(a+v_1)\cdots(a+v_g)=\mp 4$, which leaves only finitely many
possibilities for $a$ and the $v_j$. It is easy to see that the only
two examples for $g>1$ and $ 0 < v_1 < v_2 <\dots $ are
$$
  v=\{2,5\},\ a=-1 \qquad\text{and}\qquad v=\{3,4\},\ a=-2\>.
$$
(There is also $v=\{1, 3, 4\},\ a=-2$ but the resulting curve is singular.)
Finally, there is a case for $g=3$ where the last factor of \eqref{rufac}
splits into two quadratic factors over the rationals, namely $v=\{1,5,6\}$.
Thus we have found three examples for which
the 2-torsion polynomial $t(x)$ has $g+2$ rational factors,
all with a `$+$'-sign in \eqref{rueq},
$$
\begin{array}{ll}
  v=\{2,5\},   & -4t(x) = (x\!+\!1)(2x\!+\!1)(5x\!+\!1)(4x^2\!+\!6x\!+\!1)\>,\cr
  v=\{3,4\},   & -4t(x) = (2x\!+\!1)(3x\!+\!1)(4x\!+\!1)(2x^2\!+\!5x\!+\!1)\>,\cr
  v=\{1,5,6\}, & -4t(x) = (x\!+\!1)(5x\!+\!1)(6x\!+\!1)(x^2\!+\!6x\!+\!1)(4x^2\!+\!6x\!+\!1)\>.\cr
\end{array}
$$
In each of the examples, if we associate to these factors
elements $M_1,\dots,M_{g+2}$ of $\INT$ in the order that the factors are
written, we can leave out $M_{g+2}$ and still expect a non-trivial relation
between the remaining elements. Such a relation was indeed found
numerically in all three cases.

The results are summarized in Tables~\ref{tru2} (genus 2),
\ref{tru3plus}~(genus 3), \ref{tru3minus}~(genus 3) and \ref{tru4}~(genus~4).
The entries in the tables are sorted according to the conductor.
Unfortunately, we cannot compute the numerical values of $ L^*(0) $ for
$ g\geq 5 $ for these curves because the conductors get too large.

\begin{remark}\label{sequel}
The second author established in \cite{dJ} a limit formula
for the regulator of $ M_1,\dots,M_g $ (corresponding to $ v_1x+1,\dots,v_gx+1 $)
for the curves in Example~\ref{extrue} when $ v_1,\dots,v_{g-1} $ are fixed and $ |v_g| $ goes to infinity.
It follows that, for every $ g\geq2 $, there are infinitely many of such curves for which that regulator does not vanish
and $ \rk \INT \geq g $.
It is also determined precisely when the curves in Example~\ref{extrue} are isomorphic over $ \Q $ (see Remark~6.9 of \cite{dJ}),
and shown that the regulator does not vanish (and hence $ \rk \INT \geq g $) for infinitely many isomorphism classes over~$ \Qbar $.

Such results are also proven for certain hyperelliptic curves over $ \Q $ that are similar to the ones in Example~\ref{extrue}, as well
as for certain elliptic curves over arbitrary real quadratic fields.
\end{remark}

\begin{remark}
All but two of the curves in our tables have distinct conductors and are therefore pairwise non-isomorphic.
The two exceptions are the curve $ b=1 $ of Table~\ref{tab6} and the curve $ v = \{3,4\} $ of Table~\ref{tru2},
which both have conductor 816.  These curves are in fact isomorphic,
by mapping $ (x,y) $ to  $ (-\frac{x}{2x+1} , \frac{y}{(2x+1)^3} ) $.
\end{remark}

\begin{remark}\label{seven}
As mentioned in the introduction, in a number of cases we gave more than $ g $ elements of $ \INT $ and
hence by Conjecture~\ref{beilinsonconjecture} expected to find a linear relation among them.  In
all but seven cases (two in each of Tables~\ref{tab5} and \ref{tru2}, and one in each of Tables~\ref{tab6}, \ref{tab7}~and~\ref{tru3plus}),
Proposition~\ref{Vproposition} and~\eqref{urelation} of Proposition~\ref{relprop} suffice to show
that the subgroup $ \Lambda $ of $ \INT $ has at most the predicted rank $ g $.
Also in these seven cases, the regulator calculations strongly indicate that
$ \Lambda $ has rank $ g $, as was worked out in detail for one of the curves in Example~\ref{exgen2}.  
We did not try to prove the expected relation between the $ M_i $, and therefore that the rank of $ \Lambda $
is as predicted by Conjecture~\ref{beilinsonconjecture}, for any of the seven cases.
Each could, of course, be checked by a finite calculation if it is true: for instance, the conjectural
relation~\eqref{Mrelation} says that some non-zero multiple of the element
$$
 \biggl( \frac{y^2}{x^5} \biggr)  \otimes  \biggl(
  \frac{(1-x)^{41} (1-4x)^{56}} {(1-6x+11x^2-x^3)^{44}} \biggr) 
$$
in the tensor square of the multiplicative group of the function field of the curve
$y^2 - (7x^2-6x+1)y + x^5 = 0$ is a combination of tensor products
of the form  $u_i \otimes (1-u_i)$, and if this is indeed true then 
one can obviously prove it simply by exhibiting the elements $u_i$.
But in some sense it is precisely the fact that there is no visible
reason for such a relation, apart from the fact that the rank is not
expected to exceed 2 in this case, that provides the strongest experimental support
for Beilinson's conjecture.
\end{remark}

\begin{remark}
One reason we have given detailed
results for so many individual curves---apart from the fact that each
example is non-trivial to find and to calculate and that it therefore
seemed worth listing them completely---is that the experimental data about the 
values of $L^*(C,0)/R(\Lambda)$ may be useful in the future for the
formulation or numerical verification of conjectures about this rational number analogous to
the Lichtenbaum conjectures in the case of $K_n$ of number fields.
We will not make any attempt to give such conjectures.  We do note,
however, that we often find an integer for $ L^*(C,0)/R(\Lambda) $.
So it might be that $ L^*(C,0)/R(\INT) $ is always
an integer, and that the denominators that occur in $  L^*(C,0)/R(\Lambda) $
for some of our examples can be explained by the fact that 
$\Lambda$ is not the full group $\INT$ because
$ L^*(C,0)/R(\INT) = (\INT : \Lambda) \cdot L^*(C,0)/R(\Lambda) $.
This would be similar to what happens in the discussion
of the case of $ b=-1 $ 
in Example~\ref{exgen2}, where
$ \Lambda\iM = \langle M_1, M_2, M_3 \rangle $ leads to $L^*(C,0)/R(\Lambda\iM)\>\>{\buildrel?\over=}\>\>1/2$,
but the larger lattice $ \Lambda = \langle M_1, M_2, M_3, \specelt\,\rangle $,
with $ (\Lambda : \Lambda_M ) = 2 $, gives $L^*(C,0)/R(\Lambda)\>\>{\buildrel?\over=}\>\>1$.

Unfortunately, we have not been able to test this
as we do not know any systematic ways to construct
elements of $\INT$ other than the ones we discussed.
But it would be interesting to do this,
especially for the curve defined by
\[
y^2 + (5 x^3 - 13 x^2 + 7x - 1) y + x^7 = 0
\]
(which is the example with $ a=1 $ and $ b=2 $ of Table~\ref{tab7}),
where the denominator is the relatively large prime 19.
\end{remark}



\def\regquo{\vphantom{\sum_{g_{g_g}}^{\int^1}}{L^*(0)}/{R(\Lambda)} \cr}
\def\extable#1{%
\vskip-20pt
$$
\begin{array}{|r|c|c|r|c|c|}
\hline\hfill #1\hfill & \text{Conductor} & \Lambda & \hfill{L^*(0)}\hfill & \regquo \hline}
\def\exendtable{\hline\end{array}$$}
\long\def\skipeverything#1\endskip{}
\def\ph{}

\begin{table}[h]
\caption{Genus 2 curves $y^2 + ((4b-3)x^2-(5b-1)x+b) y + x^5 = 0$}
\label{tab5}
\extable{b}
\!-\!10 & 2^2\!\cdot\!3\!\cdot\!5^2\!\cdot\!7\!\cdot\!11^2\!\cdot\!1031 & \begin{array}{c}\langle M_1, M_2\rangle\end{array} & 161973.28326750 & 2^2\!\cdot\!3\!\cdot\!5^2\!\cdot\!7\cr
\!-\!9 & 3^3\!\cdot\!23\!\cdot\!8461 & \begin{array}{c}\langle M_1, M_2\rangle\end{array} & \ph\ph2647.35530780278 & 2^2\!\cdot\!3^2\cr
\!-\!8 & 2^2\!\cdot\!3\!\cdot\!61\!\cdot\!6113 & \begin{array}{c}\langle M_1, M_2\rangle\end{array} & \ph\ph2402.40716016213 & 3\!\cdot\!23/2\cr
\!-\!7 & 3\!\cdot\!7^2\!\cdot\!53\!\cdot\!4243 & \begin{array}{c}\langle M_1, M_2\rangle\end{array} & \ph16142.83542365698 & 13\!\cdot\!19\cr
\!-\!6 & 2^2\!\cdot\!3^3\!\cdot\!5\!\cdot\!2797 & \begin{array}{c}\langle M_1, M_2\rangle\end{array} & \ph\ph1090.93388314410 & 2\!\cdot\!3^2\cr
\!-\!5 & 3\!\cdot\!5^2\!\cdot\!37\!\cdot\!1721 & \begin{array}{c}\langle M_1, M_2\rangle\end{array} & \ph\ph1602.46040140666 & 29\cr
\!-\!4 & 2^2\!\cdot\!3\!\cdot\!29\!\cdot\!31^2 & \begin{array}{c}\langle M_1, M_2\rangle\end{array} & \ph\ph\ph196.40102935125 & 2^2\cr
\!-\!3 & 3^3\!\cdot\!7\!\cdot\!463 & \begin{array}{c}\langle M_1, M_2\rangle\end{array} & \ph\ph\ph\ph41.79313813775 & 1\cr
\!-\!2 & 2^2\!\cdot\!3\!\cdot\!13\!\cdot\!173 & \begin{array}{c}\langle M_1, M_2\rangle\end{array} & \ph\ph\ph\ph16.33803025573 & 1/2\cr
\!-\!1 & 3\!\cdot\!5\!\cdot\!37 & \begin{array}{c}\langle M_1, M_2, M_3, \specelt\,\rangle\cr^{_{ 41M_1 +56M_2 -44M_3 =0}}\end{array} & \ph\ph\ph\ph\ph0.22831231665 & 1\cr 
1 & 3\!\cdot\!11^2 & \begin{array}{c}\langle M_1, M_2, M_3, \specelt\,\rangle\cr^{_{ 4M_1 -26M_2 +29M_3 =0}}\end{array} & \ph\ph\ph\ph\ph0.12984963472 & 1/2\cr 
2 & 2^2\!\cdot\!3\!\cdot\!13\!\cdot\!19 & \begin{array}{c}\langle M_1, M_2\rangle\end{array} & \ph\ph\ph\ph\ph1.90319317513 & 1/2^2\!\cdot\!5\cr
3 & 3^3\!\cdot\!47 & \begin{array}{c}\langle M_1, M_2\rangle\end{array} & \ph\ph\ph\ph\ph0.62178975664 & 1/2\!\cdot\!3^3\cr
4 & 2^2\!\cdot\!3\!\cdot\!5\!\cdot\!7\!\cdot\!223 & \begin{array}{c}\langle M_1, M_2\rangle\end{array} & \ph\ph\ph\ph42.90813731246 & 1\cr
5 & 3\!\cdot\!5^2\!\cdot\!43\!\cdot\!569 & \begin{array}{c}\langle M_1, M_2\rangle\end{array} & \ph\ph\ph802.87262799199 & 2^4\cr
6 & 2^2\!\cdot\!3^3\!\cdot\!17^2\!\cdot\!67 & \begin{array}{c}\langle M_1, M_2\rangle\end{array} & \ph\ph1575.43548439889 & 2^2\!\cdot\!7\cr
7 & 3\!\cdot\!7^2\!\cdot\!59\!\cdot\!1987 & \begin{array}{c}\langle M_1, M_2\rangle\end{array} & \ph\ph6154.43484637856 & 2^2\!\cdot\!5^2\cr
8 & 2^2\!\cdot\!3\!\cdot\!67\!\cdot\!3167 & \begin{array}{c}\langle M_1, M_2\rangle\end{array} & \ph\ph1788.15229247202 & 3^3\cr
9 & 3^3\!\cdot\!5\!\cdot\!4733 & \begin{array}{c}\langle M_1, M_2\rangle\end{array} & \ph\ph\ph281.80083335457 & 2^2\cr
10 & 2^2\!\cdot\!3\!\cdot\!5^2\!\cdot\!23\!\cdot\!83\!\cdot\!293 & \begin{array}{c}\langle M_1, M_2\rangle\end{array} & \ph85596.822531781 & 2^7\!\cdot\!3^2\cr
\exendtable
\end{table}

\begin{table}[h]
\ \bigskip\ \bigskip
\caption{Genus 2 curves $y^2 + (2x^3-4bx^2-x+b) y + x^6 = 0$}
\label{tab6}
\extable{b}
1 & 2^4\!\cdot\!3\!\cdot\!17 & \begin{array}{c}\langle M_1,M_2,M_3,2\specelt\,\rangle\cr^{_{ 10M_1 -2M_2 +M_3 =0}}\end{array} & 0.35283625318 & 1/2^3\cr 
2 & 2^5\!\cdot\!3\!\cdot\!5^2\!\cdot\!13 & \begin{array}{c}\langle M_2,M_3\rangle\end{array} & 22.9767849707 & 1/2\cr
3 & 2^4\!\cdot\!3^2\!\cdot\!5^2\!\cdot\!7\!\cdot\!29 & \begin{array}{c}\langle M_2,M_3\rangle\end{array} & 347.644931188 & 2\!\cdot\!3\cr
4 & 2^5\!\cdot\!3\!\cdot\!7\!\cdot\!257 & \begin{array}{c}\langle M_2,M_3\rangle\end{array} & 134.428839855 & 2\cr
5 & 2^4\!\cdot\!3\!\cdot\!5^2\!\cdot\!11\!\cdot\!401 & \begin{array}{c}\langle M_2,M_3\rangle\end{array} & 2694.74551646 & 2^2\!\cdot\!3^2\cr
6 & 2^5\!\cdot\!3^2\!\cdot\!11\!\cdot\!13\!\cdot\!577 & \begin{array}{c}\langle M_2,M_3\rangle\end{array} & 20021.3775652 & 2\!\cdot\!3\!\cdot\!41\cr
8 & 2^5\!\cdot\!3\!\cdot\!5^2\!\cdot\!17\!\cdot\!41 & \begin{array}{c}\langle M_2,M_3\rangle\end{array} & 1106.79592308 & 2^2\!\cdot\!3\cr
10 & 2^5\!\cdot\!3\!\cdot\!5^2\!\cdot\!7\!\cdot\!19\!\cdot\!1601 & \begin{array}{c}\langle M_2,M_3\rangle\end{array} & 416121.1462 & 2^2\!\cdot\!3\!\cdot\!7^3\cr
12 & 2^5\!\cdot\!3^2\!\cdot\!5^2\!\cdot\!23\!\cdot\!461 & \begin{array}{c}\langle M_2,M_3\rangle\end{array} & 58663.7341258 & 2^2\!\cdot\!3^3\!\cdot\!5\cr
13 & 2^4\!\cdot\!3\!\cdot\!5^2\!\cdot\!13^2\!\cdot\!541 & \begin{array}{c}\langle M_2,M_3\rangle\end{array} & 46380.28556 & 2\!\cdot\!3^2\!\cdot\!23\cr
14 & 2^5\!\cdot\!3\!\cdot\!7^2\!\cdot\!29\!\cdot\!3137 & \begin{array}{c}\langle M_2,M_3\rangle\end{array} & 322837.4973 & 2\!\cdot\!3\!\cdot\!467\cr
\exendtable
\end{table}

\begin{table}[h]
\caption{Genus 3 curves
\hbox{$y^2 + ((4b\!-\!3)x^3-(4a\!+\!5b\!-\!1)x^2+(5a\!+\!b)x-a) y + x^7 = 0$}}
\label{tab7}
\extable{a,b}
\!-\!1,\!-\!6 & 3\!\cdot\!5\!\cdot\!13\!\cdot\!62773 & \begin{array}{c}\langle M_1,M_2,M_3,\specelt\,\rangle\end{array} & 157.32845991764 & 2\!\cdot\!3\cr
\!-\!1,\!-\!5 & 2\!\cdot\!3\!\cdot\!5\!\cdot\!34543 & \begin{array}{c}\langle M_1,M_2,M_3,\specelt\,\rangle\end{array} & \ph13.71462235525 & 1\cr
\!-\!1,\!-\!4 & 3^3\!\cdot\!7^2\!\cdot\!73 & \begin{array}{c}\langle M_1,M_2,M_3,\specelt\,\rangle\end{array} & \ph\ph0.99948043865 & 1/2\!\cdot\!3\cr
\!-\!1,\!-\!2 & 3\!\cdot\!19\!\cdot\!1051 & \begin{array}{c}\langle M_1,M_2,M_3,\specelt\,\rangle\end{array} & \ph\ph0.68841617094 & 1/2\!\cdot\!7\cr
\!-\!1,0 & 3\!\cdot\!5\!\cdot\!7\!\cdot\!997 & \begin{array}{c}\langle M_1,M_2,M_3,\specelt\,\rangle\end{array} & \ph\ph1.30028297708 & 1/2\!\cdot\!3\cr
\!-\!1,1 & 2\!\cdot\!3\!\cdot\!43\!\cdot\!599 & \begin{array}{c}\langle M_1,M_2,M_3,\specelt\,\rangle\end{array} & \ph\ph1.85769147359 & 1/2\cr
\!-\!1,2 & 3^3\!\cdot\!17\!\cdot\!199 & \begin{array}{c}\langle M_1,M_2,M_3,\specelt\,\rangle\end{array} & \ph\ph0.95009772024 & 1/2\!\cdot\!3\!\cdot\!5\cr
\!-\!1,3 & 2\!\cdot\!3\!\cdot\!7^2\!\cdot\!59\!\cdot\!599 & \begin{array}{c}\langle M_1,M_2,M_3,\specelt\,\rangle\end{array} & 117.81139857650 & 2^2\cr
1,0 & 3\!\cdot\!29^2\!\cdot\!71 & \begin{array}{c}\langle M_1,M_2,M_3,\specelt\,\rangle\end{array} & \ph\ph1.91931648711 & 1/3\cr
1,2 & 3\!\cdot\!13\!\cdot\!971 & \begin{array}{c}\langle M_1,M_2,M_3,\specelt\,\rangle\end{array} & \ph\ph0.41650834991 & 1/19\cr
1,3 & 2\!\cdot\!3\!\cdot\!5^2\!\cdot\!229 & \begin{array}{c}\langle M_1,M_2,M_{3,1},M_{3,2},\specelt\,\rangle\cr^{_{ 36M_1 -111M_2 -41M_{3,1} +71M_{3,2} =0}}\end{array} & \ph\ph0.33653518886 & 1\cr
1,4 & 3^2\!\cdot\!7877 & \begin{array}{c}\langle M_1,M_2,M_3,\specelt\,\rangle\end{array} & \ph\ph0.84449758754 & 2/3\!\cdot\!5\cr
1,5 & 2\!\cdot\!3\!\cdot\!11\!\cdot\!44071 & \begin{array}{c}\langle M_1,M_2,M_3,\specelt\,\rangle\end{array} & \ph39.13475965836 & 2\cr
1,6 & 3\!\cdot\!5^2\!\cdot\!19\!\cdot\!46273 & \begin{array}{c}\langle M_1,M_2,M_3,\specelt\,\rangle\end{array} & 631.40978672221 & 2^2\!\cdot\!5\cr
1,7 & 2\!\cdot\!3^2\!\cdot\!479\!\cdot\!2011 & \begin{array}{c}\langle M_1,M_2,M_3,\specelt\,\rangle\end{array} & 170.28405530697 & 2^2\cr
\exendtable
\end{table}

\begin{table}[h]
\ \bigskip\ \bigskip
\caption{Genus 4 curves $y^2 + f(x) y + x^{10} = 0$}
\label{tab10}
\vskip-20pt
$$
\begin{array}{|l|c|r|c|c|}
\hline
\hfill f(x)\hfill & \text{Conductor} & \hfill{L^*(0)}\hfill & \regquo
\hline
2x^5\!+\!2x^4\!+\!x^3\!+\!x^2\!-\!3x\!-\!1 & 2^{11}\!\cdot\!5^3\!\cdot\!19\!\cdot\!29 & \ph35.85879769 & 1/2\cr
2x^5\!+\!2x^4\!+\!2x^3\!-\!3x^2\!-\!2x\!+\!1 & 2^{11}\!\cdot\!3^2\!\cdot\!17\!\cdot\!59 & \ph\ph5.336928011 & 1/2\cr
2x^5\!+\!3x^4\!-\!3x^3\!-\!2x^2\!-\!x\!-\!1 & 2^3\!\cdot\!3^3\!\cdot\!5\!\cdot\!19\!\cdot\!331 & \ph\ph1.865694255 & 1/2^2\!\cdot\!3\cr
2x^5\!+\!3x^4\!-\!x^3\!-\!4x^2\!-\!3x\!+\!1 & 2^3\!\cdot\!3^4\!\cdot\!5\!\cdot\!7\!\cdot\!20759 & 126.4283012 & 1\cr
2x^5\!+\!3x^4\!+\!x^3\!-\!3x\!-\!1 & 2^4\!\cdot\!3^3\!\cdot\!7\!\cdot\!11\!\cdot\!4793 & \ph41.29358643 & 1/2\cr
2x^5\!+\!4x^4\!-\!3x^3\!-\!2x\!+\!1 & 2^4\!\cdot\!5^3\!\cdot\!7\!\cdot\!11\!\cdot\!103 & \ph\ph3.546483598 & 1/2^3\cr
2x^5\!+\!4x^4\!-\!x^3\!-\!3x^2\!-\!3x\!-\!1 & 2^{10}\!\cdot\!3\!\cdot\!7\!\cdot\!1051 & \ph\ph6.484247251 & 1/2^3\cr
2x^5\!+\!4x^4\!-\!x^3\!-\!3x^2\!-\!x\!+\!1 & 2^{12}\!\cdot\!3^6\!\cdot\!13 & \ph11.50901911 & 1/2^2\cr
2x^5\!+\!4x^4\!+\!3x^3\!-\!5x^2\!-\!5x\!-\!1 & 2^{12}\!\cdot\!3\!\cdot\!5\!\cdot\!19\!\cdot\!79 & \ph27.69939565 & 1/2\cr
2x^5\!+\!4x^4\!+\!3x^3\!+\!3x^2\!-\!x\!-\!1 & 2^{12}\!\cdot\!3\!\cdot\!5^2\!\cdot\!23\!\cdot\!43 & \ph89.28895569 & 1\cr
2x^5\!+\!4x^4\!+\!5x^3\!+\!2x^2\!+\!2x\!+\!1 & 2^4\!\cdot\!3^3\!\cdot\!7^2\!\cdot\!379 & \ph\ph2.157167657 & 1/2^4\cr
\exendtable
\end{table}

\begin{table}[h]
\ \bigskip\ \bigskip
\caption{Genus 5 curves $y^2 + f(x) y + x^{12} = 0$}
\label{tab12}
\vskip-20pt
$$
\begin{array}{|l|c|r|c|c|}
\hline
\hfill f(x) \hfill & \text{Conductor} & \hfill{L^*(0)}\hfill & \regquo
\hline
2x^6 \!+\! 2x^5 \!-\! 4x^4 \!-\! 3x^3 \!-\! 2x^2 \!+\! 4x \!-\! 1 & 2^8\!\cdot\!3^4\!\cdot\!5\!\cdot\!7\!\cdot\!11\!\cdot\!19 & \ph0.97906446422637 & 1/2^2\!\cdot\!3^3\cr
2x^6 \!+\! 4x^5 \!-\! 5x^4 \!-\! 3x^3 \!+\! 2x^2 \!-\! x \!-\! 1 & 2^7\!\cdot\!3^2\!\cdot\!5^2\!\cdot\!107\!\cdot\!139 & \ph2.86707608488323 & 1/2^2\!\cdot\!3\cr
2x^6 \!+\! 6x^5 \!-\! 5x^3 \!-\! 3x^2 \!+\! x \!+\! 1 & 2^{16}\!\cdot\!3^2\!\cdot\!5\!\cdot\!13^2 & \ph3.21518014215484 & 1/2^2\!\cdot\!3\cr
2x^6 \!+\! 2x^5 \!+\! 2x^4 \!-\! x^3 \!-\! 3x^2 \!-\! 3x \!-\! 1 & 2^{16}\!\cdot\!3\!\cdot\!5^3\!\cdot\!31 & \ph4.04748393920751 & 1/2^3\!\cdot\!3\cr
2x^6 \!+\! 2x^5 \!+\! x^3 \!-\! 3x^2 \!-\! x \!+\! 1 & 2^{16}\!\cdot\!3^4\!\cdot\!5\!\cdot\!181 & 28.41118880946 & 1/3\cr
\hline
\end{array}
$$
\end{table}

\begin{table}[h]
\caption{Genus 2 curves $y^2+(2x^3\!+\!(v_1x\!+\!1)(v_2x\!+\!1))y+x^6=0\>$}
\label{tru2}
\extable{v}
2, 5 & 2^3\!\cdot\!3\!\cdot\!5^2 & \begin{array}{c}\langle M_1, M_2, M_3,2\specelt\,\rangle\cr ^{_{39M_1-26M_2+M_3=0}} \end{array} & 0.246481356638 & 1/2^2\cr 
2, 4 & 2^6\!\cdot\!11 & \begin{array}{c}\langle M_1, M_2,2\specelt\,\rangle \end{array} & 0.339189596082 & 1/2^4\!\cdot\!3\cr 
3, 4 & 2^4\!\cdot\!3\!\cdot\!17 & \begin{array}{c}\langle M_1, M_2, M_3,2\specelt\,\rangle\cr ^{_{9M_1-10M_2+2M_3=0}} \end{array} & 0.352836253176 & 1/2^3\cr 
1, 8 & 2\!\cdot\!7\!\cdot\!59 & \begin{array}{c}\langle M_1, M_2,2\specelt\,\rangle \end{array} & 0.412953201772 & 1/2^3\!\cdot\!7\cr 
1, 2 & 2^3\!\cdot\!107 & \begin{array}{c}\langle M_1, M_2,2\specelt\,\rangle \end{array} & 0.372337494391 & 1/2^3\cr 
1, 4 & 2^3\!\cdot\!3\!\cdot\!53 & \begin{array}{c}\langle M_1, M_2,2\specelt\,\rangle \end{array} & 0.609341856994 & 1/2^4\cr 
8, 9 & 2\!\cdot\!3\!\cdot\!223 & \begin{array}{c}\langle M_1, M_2,2\specelt\,\rangle \end{array} & 0.751762610683 & 1/2^2\!\cdot\!7\cr 
2, 3 & 2^3\!\cdot\!3\!\cdot\!59 & \begin{array}{c}\langle M_1, M_2,2\specelt\,\rangle \end{array} & 0.682403286269 & 1/2^2\!\cdot\!3\cr 
4, 5 & 2^3\!\cdot\!5\!\cdot\!79 & \begin{array}{c}\langle M_1, M_2,2\specelt\,\rangle \end{array} & 1.466505046768 & 1/2^3\cr 
1, 3 & 2^3\!\cdot\!3\!\cdot\!139 & \begin{array}{c}\langle M_1, M_2,2\specelt\,\rangle \end{array} & 1.573797007363 & 1/2^2\cr 
1, 5 & 2^4\!\cdot\!5\!\cdot\!83 & \begin{array}{c}\langle M_1, M_2,2\specelt\,\rangle \end{array} & 3.296439983269 & 1/2^2\cr 
2, 6 & 2^6\!\cdot\!3\!\cdot\!37 & \begin{array}{c}\langle M_1, M_2,2\specelt\,\rangle \end{array} & 4.318055838622 & 1/2^2\cr 
5, 8 & 2^2\!\cdot\!3\!\cdot\!5\!\cdot\!127 & \begin{array}{c}\langle M_1, M_2,2\specelt\,\rangle \end{array} & 4.980029782458 & 1/2\!\cdot\!3\cr 
1, 9 & 2^5\!\cdot\!3\!\cdot\!149 & \begin{array}{c}\langle M_1, M_2,2\specelt\,\rangle \end{array} & 6.585355387353 & 1/2^2\cr 
3, 5 & 2^3\!\cdot\!3\!\cdot\!5\!\cdot\!229 & \begin{array}{c}\langle M_1, M_2,2\specelt\,\rangle \end{array} & 14.988121207002 & 1\cr 
3, 11 & 2^2\!\cdot\!3\!\cdot\!11\!\cdot\!229 & \begin{array}{c}\langle M_1, M_2,2\specelt\,\rangle \end{array} & 19.618673668593 & 1/2\cr 
4, 8 & 2^5\!\cdot\!997 & \begin{array}{c}\langle M_1, M_2,2\specelt\,\rangle \end{array} & 15.361499842059 & 1/2\cr 
4, 6 & 2^6\!\cdot\!3\!\cdot\!197 & \begin{array}{c}\langle M_1, M_2,2\specelt\,\rangle \end{array} & 19.481730347124 & 1\cr 
3, 6 & 2^3\!\cdot\!3^5\!\cdot\!23 & \begin{array}{c}\langle M_1, M_2,2\specelt\,\rangle \end{array} & 20.448409445891 & 1\cr 
2, 8 & 2^6\!\cdot\!3\!\cdot\!269 & \begin{array}{c}\langle M_1, M_2,2\specelt\,\rangle \end{array} & 24.399036028269 & 1\cr 
2, 10 & 2^6\!\cdot\!5^2\!\cdot\!37 & \begin{array}{c}\langle M_1, M_2,2\specelt\,\rangle \end{array} & 30.231463268350 & 1\cr 
1, 6 & 2^3\!\cdot\!3\!\cdot\!5\!\cdot\!499 & \begin{array}{c}\langle M_1, M_2,2\specelt\,\rangle \end{array} & 24.839835330300 & 3/2\cr 
6, 9 & 2^3\!\cdot\!3^5\!\cdot\!31 & \begin{array}{c}\langle M_1, M_2,2\specelt\,\rangle \end{array} & 33.600131091460 & 1\cr 
3, 7 & 2^3\!\cdot\!3\!\cdot\!7\!\cdot\!359 & \begin{array}{c}\langle M_1, M_2,2\specelt\,\rangle \end{array} & 25.126558038878 & 1\cr 
7, 8 & 2^5\!\cdot\!7^2\!\cdot\!41 & \begin{array}{c}\langle M_1, M_2,2\specelt\,\rangle \end{array} & 28.545427660372 & 3/2\cr 
5, 6 & 2^3\!\cdot\!3\!\cdot\!5\!\cdot\!733 & \begin{array}{c}\langle M_1, M_2,2\specelt\,\rangle \end{array} & 43.322633518831 & 3\cr 
3, 12 & 2^4\!\cdot\!3^5\!\cdot\!23 & \begin{array}{c}\langle M_1, M_2,2\specelt\,\rangle \end{array} & 42.043168670629 & 1\cr 
1, 7 & 2^3\!\cdot\!3\!\cdot\!7^2\!\cdot\!101 & \begin{array}{c}\langle M_1, M_2,2\specelt\,\rangle \end{array} & 59.608745714348 & 3\cr 
2, 7 & 2^3\!\cdot\!5\!\cdot\!7^2\!\cdot\!67 & \begin{array}{c}\langle M_1, M_2,2\specelt\,\rangle \end{array} & 63.124534380486 & 3\cr 
4, 7 & 2^4\!\cdot\!3\!\cdot\!7^2\!\cdot\!67 & \begin{array}{c}\langle M_1, M_2,2\specelt\,\rangle \end{array} & 76.831677676717 & 3\cr 
6, 8 & 2^6\!\cdot\!3\!\cdot\!829 & \begin{array}{c}\langle M_1, M_2,2\specelt\,\rangle \end{array} & 79.443510393961 & 3\cr 
1, 10 & 2^3\!\cdot\!3\!\cdot\!5\!\cdot\!1427 & \begin{array}{c}\langle M_1, M_2,2\specelt\,\rangle \end{array} & 88.599127797882 & 3\cr 
3, 8 & 2^5\!\cdot\!3\!\cdot\!5^2\!\cdot\!73 & \begin{array}{c}\langle M_1, M_2,2\specelt\,\rangle \end{array} & 87.679227492015 & 3\cr 
4, 9 & 2^3\!\cdot\!3\!\cdot\!5\!\cdot\!1907 & \begin{array}{c}\langle M_1, M_2,2\specelt\,\rangle \end{array} & 105.425224481666 & 3\cr 
6, 7 & 2^3\!\cdot\!3\!\cdot\!7\!\cdot\!1373 & \begin{array}{c}\langle M_1, M_2,2\specelt\,\rangle \end{array} & 101.061871259886 & 2\!\cdot\!3\cr 
2, 9 & 2^3\!\cdot\!3\!\cdot\!7\!\cdot\!11\!\cdot\!191 & \begin{array}{c}\langle M_1, M_2,2\specelt\,\rangle \end{array} & 164.650770360786 & 2\!\cdot\!3\cr 
3, 9 & 2^3\!\cdot\!3^5\!\cdot\!199 & \begin{array}{c}\langle M_1, M_2,2\specelt\,\rangle \end{array} & 197.318826602911 & 2\!\cdot\!3\cr 
5, 7 & 2^3\!\cdot\!5^2\!\cdot\!7\!\cdot\!353 & \begin{array}{c}\langle M_1, M_2,2\specelt\,\rangle \end{array} & 208.897958121197 & 3^2\cr 
5, 9 & 2^4\!\cdot\!3\!\cdot\!5\!\cdot\!2089 & \begin{array}{c}\langle M_1, M_2,2\specelt\,\rangle \end{array} & 212.013278144012 & 2\!\cdot\!3\cr 
\exendtable
\end{table}

\begin{table}[h]
\vskip20 pt
$$
\begin{array}{|r|c|c|r|c|c|}
\hline
9, 10 & 2^3\!\cdot\!3\!\cdot\!5\!\cdot\!5261 & \begin{array}{c}\langle M_1, M_2,2\specelt\,\rangle \end{array} & 275.183140339336 & 2^2\!\cdot\!3\cr 
8, 10 & 2^6\!\cdot\!5\!\cdot\!2221 & \begin{array}{c}\langle M_1, M_2,2\specelt\,\rangle \end{array} & 417.499914032374 & 13\cr 
1, 11 & 2^3\!\cdot\!5\!\cdot\!11\!\cdot\!1619 & \begin{array}{c}\langle M_1, M_2,2\specelt\,\rangle \end{array} & 294.450550333465 & 3^2\cr 
4, 12 & 2^5\!\cdot\!3\!\cdot\!83\!\cdot\!103 & \begin{array}{c}\langle M_1, M_2,2\specelt\,\rangle \end{array} & 412.225445696365 & 3^2\cr 
7, 9 & 2^3\!\cdot\!3\!\cdot\!7^2\!\cdot\!11\!\cdot\!73 & \begin{array}{c}\langle M_1, M_2,2\specelt\,\rangle \end{array} & 441.365714448209 & 3\!\cdot\!5\cr 
2, 11 & 2^3\!\cdot\!3\!\cdot\!11\!\cdot\!6053 & \begin{array}{c}\langle M_1, M_2,2\specelt\,\rangle \end{array} & 689.050402379536 & 3\!\cdot\!7\cr 
5, 12 & 2^3\!\cdot\!3\!\cdot\!5^2\!\cdot\!7^2\!\cdot\!61 & \begin{array}{c}\langle M_1, M_2,2\specelt\,\rangle \end{array} & 860.242245920864 & 2\!\cdot\!3^2\cr 
2, 12 & 2^6\!\cdot\!3\!\cdot\!5\!\cdot\!2341 & \begin{array}{c}\langle M_1, M_2,2\specelt\,\rangle \end{array} & 1267.775571888346 & 2^2\!\cdot\!3^2\cr 
7, 15 & 2\!\cdot\!3\!\cdot\!5^2\!\cdot\!7^2\!\cdot\!313 & \begin{array}{c}\langle M_1, M_2,2\specelt\,\rangle \end{array} & 1075.913430726619 & 2\!\cdot\!3^2\cr 
8, 17 & 2\!\cdot\!3\!\cdot\!17\!\cdot\!23321 & \begin{array}{c}\langle M_1, M_2,2\specelt\,\rangle \end{array} & 1201.189781491657 & 2\!\cdot\!3^2\cr 
5, 10 & 2^3\!\cdot\!5^2\!\cdot\!59\!\cdot\!263 & \begin{array}{c}\langle M_1, M_2,2\specelt\,\rangle \end{array} & 1440.017545464245 & 2^2\!\cdot\!3^2\cr 
4, 10 & 2^6\!\cdot\!3\!\cdot\!5\!\cdot\!17\!\cdot\!197 & \begin{array}{c}\langle M_1, M_2,2\specelt\,\rangle \end{array} & 1562.214360456082 & 2^3\!\cdot\!5\cr 
6, 10 & 2^6\!\cdot\!3\!\cdot\!5\!\cdot\!3797 & \begin{array}{c}\langle M_1, M_2,2\specelt\,\rangle \end{array} & 1888.148590577887 & 2^4\!\cdot\!3\cr 
4, 13 & 2^3\!\cdot\!3\!\cdot\!11\!\cdot\!13^2\!\cdot\!89 & \begin{array}{c}\langle M_1, M_2,2\specelt\,\rangle \end{array} & 2048.440414253531 & 2\!\cdot\!3\!\cdot\!7\cr 
8, 16 & 2^6\!\cdot\!109\!\cdot\!601 & \begin{array}{c}\langle M_1, M_2,2\specelt\,\rangle \end{array} & 2029.415672068518 & 2^5\cr 
4, 14 & 2^6\!\cdot\!5\!\cdot\!7^2\!\cdot\!373 & \begin{array}{c}\langle M_1, M_2,2\specelt\,\rangle \end{array} & 3298.078842289480 & 2^6\cr 
6, 12 & 2^6\!\cdot\!3^5\!\cdot\!431 & \begin{array}{c}\langle M_1, M_2,2\specelt\,\rangle \end{array} & 3779.698402282094 & 2\!\cdot\!3\!\cdot\!13\cr 
4, 11 & 2^4\!\cdot\!7\!\cdot\!11\!\cdot\!23\!\cdot\!239 & \begin{array}{c}\langle M_1, M_2,2\specelt\,\rangle \end{array} & 3065.942241509888 & 2^3\!\cdot\!3^2\cr 
3, 10 & 2^3\!\cdot\!3\!\cdot\!5\!\cdot\!7^2\!\cdot\!1307 & \begin{array}{c}\langle M_1, M_2,2\specelt\,\rangle \end{array} & 3909.737148728376 & 2^2\!\cdot\!3^3\cr 
7, 10 & 2^3\!\cdot\!3\!\cdot\!5^2\!\cdot\!7\!\cdot\!43\!\cdot\!59 & \begin{array}{c}\langle M_1, M_2,2\specelt\,\rangle \end{array} & 4653.068964642799 & 2\!\cdot\!3^2\!\cdot\!7\cr 
5, 13 & 2^5\!\cdot\!5\!\cdot\!13\!\cdot\!31\!\cdot\!263 & \begin{array}{c}\langle M_1, M_2,2\specelt\,\rangle \end{array} & 8592.680376054557 & 2^3\!\cdot\!3\!\cdot\!7\cr 
7, 16 & 2^5\!\cdot\!3\!\cdot\!7^2\!\cdot\!4493 & \begin{array}{c}\langle M_1, M_2,2\specelt\,\rangle \end{array} & 11344.672982180578 & 2^2\!\cdot\!3^2\!\cdot\!5\cr 
5, 11 & 2^3\!\cdot\!3\!\cdot\!5\!\cdot\!11\!\cdot\!26573 & \begin{array}{c}\langle M_1, M_2,2\specelt\,\rangle \end{array} & 16672.060582310110 & 2\!\cdot\!3^3\!\cdot\!7\cr 
6, 14 & 2^6\!\cdot\!3\!\cdot\!7\!\cdot\!41\!\cdot\!677 & \begin{array}{c}\langle M_1, M_2,2\specelt\,\rangle \end{array} & 20507.523505991035 & 2^4\!\cdot\!23\cr 
7, 14 & 2^3\!\cdot\!7^2\!\cdot\!117541 & \begin{array}{c}\langle M_1, M_2,2\specelt\,\rangle \end{array} & 20242.092778779142 & 2^3\!\cdot\!3^2\!\cdot\!5\cr 
3, 13 & 2^3\!\cdot\!3\!\cdot\!5^2\!\cdot\!13\!\cdot\!6553 & \begin{array}{c}\langle M_1, M_2,2\specelt\,\rangle \end{array} & 20094.811949554845 & 2\!\cdot\!3^2\!\cdot\!5^2\cr 
6, 16 & 2^6\!\cdot\!3\!\cdot\!5\!\cdot\!73\!\cdot\!773 & \begin{array}{c}\langle M_1, M_2,2\specelt\,\rangle \end{array} & 26721.552091088677 & 2^4\!\cdot\!3^3\cr 
6, 15 & 2^3\!\cdot\!3^5\!\cdot\!5\!\cdot\!59\!\cdot\!101 & \begin{array}{c}\langle M_1, M_2,2\specelt\,\rangle \end{array} & 28628.778269858471 & 2\!\cdot\!3^5\cr 
5, 14 & 2^3\!\cdot\!3\!\cdot\!5\!\cdot\!7\!\cdot\!95621 & \begin{array}{c}\langle M_1, M_2,2\specelt\,\rangle \end{array} & 32866.799798423930 & 2\!\cdot\!3\!\cdot\!101\cr 
5, 15 & 2^3\!\cdot\!3\!\cdot\!5^2\!\cdot\!29\!\cdot\!4673 & \begin{array}{c}\langle M_1, M_2,2\specelt\,\rangle \end{array} & 41114.005206032957 & 2^4\!\cdot\!3^2\!\cdot\!5\cr 
8, 18 & 2^6\!\cdot\!3\!\cdot\!5^2\!\cdot\!73\!\cdot\!353 & \begin{array}{c}\langle M_1, M_2,2\specelt\,\rangle \end{array} & 59468.612113790229 & 2^2\!\cdot\!3\!\cdot\!71\cr 
6, 13 & 2^3\!\cdot\!3\!\cdot\!7^2\!\cdot\!13\!\cdot\!17\!\cdot\!619 & \begin{array}{c}\langle M_1, M_2,2\specelt\,\rangle \end{array} & 81849.542063836685 & 2\!\cdot\!3^3\!\cdot\!29\cr 
\exendtable
\end{table}

\begin{table}[h]
\caption{Genus 3 curves $y^2+(2x^4\!+\!(v_1x\!+\!1)(v_2x\!+\!1)(v_3x\!+\!1))y+x^8=0\>$}
\label{tru3plus}
\extable{v}
2, 4, 6 & 2^8\!\cdot\!3^3\!\cdot\!13 & \begin{array}{c}\langle M_1, M_2, M_3,2\specelt\,\rangle \end{array} & 1.237888702787 & 1/2^7\cr 
1, 3, 4 & 2^8\!\cdot\!3^2\!\cdot\!7^2 & \begin{array}{c}\langle M_1, M_2, M_3,2\specelt\,\rangle \end{array} & 1.319632692018 & 1/2^6\cr 
1, 2, 4 & 2^6\!\cdot\!3^3\!\cdot\!67 & \begin{array}{c}\langle M_1, M_2, M_3,2\specelt\,\rangle \end{array} & 1.326260789773 & 1/2^5\cr 
1, 2, 3 & 2^6\!\cdot\!3^3\!\cdot\!73 & \begin{array}{c}\langle M_1, M_2, M_3,2\specelt\,\rangle \end{array} & 1.360553792128 & 1/2^4\cr 
1, 5, 6 & 2^8\!\cdot\!3^3\!\cdot\!5^3 & \begin{array}{c}\langle M_1, M_2, M_3, M_4,2\specelt\,\rangle\cr ^{_{M_1+M_2+11M_3-12M_4=0}} \end{array} & 9.556528296211 & 1/2\cr 
1, 3, 5 & 2^6\!\cdot\!3^3\!\cdot\!5\!\cdot\!179 & \begin{array}{c}\langle M_1, M_2, M_3,2\specelt\,\rangle \end{array} & 15.763316182605 & 1/2^3\cr 
3, 4, 7 & 2^8\!\cdot\!3^2\!\cdot\!7\!\cdot\!113 & \begin{array}{c}\langle M_1, M_2, M_3,2\specelt\,\rangle \end{array} & 21.282782316338 & 1/2^3\cr 
1, 2, 5 & 2^8\!\cdot\!3^3\!\cdot\!5\!\cdot\!53 & \begin{array}{c}\langle M_1, M_2, M_3,2\specelt\,\rangle \end{array} & 18.574599809025 & 1/2^2\cr 
2, 6, 8 & 2^{14}\!\cdot\!3^2\!\cdot\!17 & \begin{array}{c}\langle M_1, M_2, M_3,2\specelt\,\rangle \end{array} & 41.533490724724 & 1/2^3\cr 
2, 3, 5 & 2^6\!\cdot\!3^2\!\cdot\!5^3\!\cdot\!89 & \begin{array}{c}\langle M_1, M_2, M_3,2\specelt\,\rangle \end{array} & 72.152711581916 & 1\cr 
1, 2, 6 & 2^6\!\cdot\!3^3\!\cdot\!5^3\!\cdot\!43 & \begin{array}{c}\langle M_1, M_2, M_3,2\specelt\,\rangle \end{array} & 105.836805455409 & 1\cr 
1, 3, 6 & 2^6\!\cdot\!3^2\!\cdot\!5\!\cdot\!47\!\cdot\!73 & \begin{array}{c}\langle M_1, M_2, M_3,2\specelt\,\rangle \end{array} & 98.432671578138 & 1/2\cr 
2, 3, 6 & 2^6\!\cdot\!3^2\!\cdot\!7^2\!\cdot\!373 & \begin{array}{c}\langle M_1, M_2, M_3,2\specelt\,\rangle \end{array} & 113.826557165808 & 1\cr 
1, 3, 7 & 2^6\!\cdot\!3^2\!\cdot\!7\!\cdot\!2837 & \begin{array}{c}\langle M_1, M_2, M_3,2\specelt\,\rangle \end{array} & 137.788942178887 & 1\cr 
1, 4, 6 & 2^6\!\cdot\!3^2\!\cdot\!5\!\cdot\!5669 & \begin{array}{c}\langle M_1, M_2, M_3,2\specelt\,\rangle \end{array} & 198.333924314685 & 3/2\cr 
1, 2, 9 & 2^7\!\cdot\!3^3\!\cdot\!7\!\cdot\!883 & \begin{array}{c}\langle M_1, M_2, M_3,2\specelt\,\rangle \end{array} & 206.791570180633 & 1\cr 
3, 4, 8 & 2^5\!\cdot\!3^3\!\cdot\!5\!\cdot\!7^2\!\cdot\!101 & \begin{array}{c}\langle M_1, M_2, M_3,2\specelt\,\rangle \end{array} & 219.415523333470 & 1\cr 
4, 5, 9 & 2^8\!\cdot\!3^3\!\cdot\!5^3\!\cdot\!29 & \begin{array}{c}\langle M_1, M_2, M_3,2\specelt\,\rangle \end{array} & 278.222780993318 & 1\cr 
2, 4, 10 & 2^{11}\!\cdot\!3^3\!\cdot\!5\!\cdot\!101 & \begin{array}{c}\langle M_1, M_2, M_3,2\specelt\,\rangle \end{array} & 375.379161057633 & 1\cr 
1, 2, 8 & 2^6\!\cdot\!3^3\!\cdot\!7\!\cdot\!29\!\cdot\!103 & \begin{array}{c}\langle M_1, M_2, M_3,2\specelt\,\rangle \end{array} & 429.330911878658 & 5/2\cr 
1, 4, 9 & 2^7\!\cdot\!3^2\!\cdot\!5\!\cdot\!9281 & \begin{array}{c}\langle M_1, M_2, M_3,2\specelt\,\rangle \end{array} & 554.543543287698 & 2\cr 
1, 3, 9 & 2^6\!\cdot\!3^2\!\cdot\!137437 & \begin{array}{c}\langle M_1, M_2, M_3,2\specelt\,\rangle \end{array} & 831.548080326194 & 2^2\cr 
1, 5, 8 & 2^7\!\cdot\!3^3\!\cdot\!5\!\cdot\!7\!\cdot\!659 & \begin{array}{c}\langle M_1, M_2, M_3,2\specelt\,\rangle \end{array} & 798.968120685204 & 3\cr 
1, 4, 7 & 2^8\!\cdot\!3^2\!\cdot\!7\!\cdot\!5879 & \begin{array}{c}\langle M_1, M_2, M_3,2\specelt\,\rangle \end{array} & 930.405544167007 & 5\cr 
1, 2, 7 & 2^6\!\cdot\!3^3\!\cdot\!5\!\cdot\!7\!\cdot\!1999 & \begin{array}{c}\langle M_1, M_2, M_3,2\specelt\,\rangle \end{array} & 1105.646745724194 & 2^3\cr 
1, 3, 8 & 2^7\!\cdot\!3^3\!\cdot\!5^3\!\cdot\!7\!\cdot\!41 & \begin{array}{c}\langle M_1, M_2, M_3,2\specelt\,\rangle \end{array} & 1395.300068143365 & 2^3\cr 
1, 4, 8 & 2^5\!\cdot\!3^3\!\cdot\!7\!\cdot\!21773 & \begin{array}{c}\langle M_1, M_2, M_3,2\specelt\,\rangle \end{array} & 1402.036187242234 & 2\!\cdot\!3\cr 
2, 4, 8 & 2^{11}\!\cdot\!3^3\!\cdot\!4093 & \begin{array}{c}\langle M_1, M_2, M_3,2\specelt\,\rangle \end{array} & 3341.033608754311 & 2^2\!\cdot\!3\cr 
\exendtable
\end{table}

\begin{table}[h]
\caption{Genus 3 curves $y^2+(2x^4\!-\!(v_1x\!+\!1)(v_2x\!+\!1)(v_3x\!+\!1))y+x^8=0\>$}
\label{tru3minus}
\extable{v}
1, 2, 3 & 2^6\!\cdot\!3\!\cdot\!5^3\!\cdot\!11 & \begin{array}{c}\langle M_1, M_2, M_3, 2\specelt\,\rangle \end{array} & 2.288179184223 & 1/2^3\cr 
1, 2, 4 & 2^6\!\cdot\!3\!\cdot\!5^3\!\cdot\!23 & \begin{array}{c}\langle M_1, M_2, M_3, 2\specelt\,\rangle \end{array} & 5.368903572985 & 1/2^3\cr 
1, 2, 5 & 2^7\!\cdot\!3\!\cdot\!5\!\cdot\!433 & \begin{array}{c}\langle M_1, M_2, M_3, 2\specelt\,\rangle \end{array} & 8.855164567689 & 1/2^3\cr 
2, 4, 6 & 2^{12}\!\cdot\!3\!\cdot\!5^3 & \begin{array}{c}\langle M_1, M_2, M_3, 2\specelt\,\rangle \end{array} & 19.749659522297 & 1/2^3\cr 
3, 4, 8 & 2^5\!\cdot\!3\!\cdot\!5\!\cdot\!7\!\cdot\!19\!\cdot\!31 & \begin{array}{c}\langle M_1, M_2, M_3, 2\specelt\,\rangle \end{array} & 21.559124747863 & 1/2^4\cr 
1, 4, 5 & 2^{10}\!\cdot\!3\!\cdot\!5\!\cdot\!7\!\cdot\!23 & \begin{array}{c}\langle M_1, M_2, M_3, 2\specelt\,\rangle \end{array} & 22.811436896592 & 1/2^2\cr 
1, 3, 5 & 2^6\!\cdot\!3\!\cdot\!5\!\cdot\!3797 & \begin{array}{c}\langle M_1, M_2, M_3, 2\specelt\,\rangle \end{array} & 42.967686126549 & 1/2\cr 
2, 3, 7 & 2^7\!\cdot\!3\!\cdot\!5\!\cdot\!7\!\cdot\!397 & \begin{array}{c}\langle M_1, M_2, M_3, 2\specelt\,\rangle \end{array} & 53.453545184856 & 1/2^2\cr 
1, 2, 10 & 2^6\!\cdot\!3\!\cdot\!5\!\cdot\!11\!\cdot\!523 & \begin{array}{c}\langle M_1, M_2, M_3, 2\specelt\,\rangle \end{array} & 57.955617125753 & 1/2^2\cr 
2, 3, 6 & 2^6\!\cdot\!3^2\!\cdot\!5^3\!\cdot\!7\!\cdot\!17 & \begin{array}{c}\langle M_1, M_2, M_3, 2\specelt\,\rangle \end{array} & 114.624730232758 & 3/2^2\cr 
1, 2, 8 & 2^6\!\cdot\!3\!\cdot\!5^3\!\cdot\!7\!\cdot\!67 & \begin{array}{c}\langle M_1, M_2, M_3, 2\specelt\,\rangle \end{array} & 122.595870053997 & 3/2^2\cr 
1, 2, 6 & 2^6\!\cdot\!3\!\cdot\!5\!\cdot\!14731 & \begin{array}{c}\langle M_1, M_2, M_3, 2\specelt\,\rangle \end{array} & 125.574116230804 & 5/2^2\cr 
1, 5, 6 & 2^7\!\cdot\!3\!\cdot\!5^2\!\cdot\!41^2 & \begin{array}{c}\langle M_1, M_2, M_3, 2\specelt\,\rangle \end{array} & 140.322800930811 & 1\cr 
1, 2, 9 & 2^6\!\cdot\!3\!\cdot\!5^3\!\cdot\!7\!\cdot\!113 & \begin{array}{c}\langle M_1, M_2, M_3, 2\specelt\,\rangle \end{array} & 196.848432921070 & 1\cr 
2, 3, 5 & 2^6\!\cdot\!3^3\!\cdot\!5\!\cdot\!13^3 & \begin{array}{c}\langle M_1, M_2, M_3, 2\specelt\,\rangle \end{array} & 190.825343100486 & 2\cr 
3, 4, 9 & 2^7\!\cdot\!3^2\!\cdot\!5\!\cdot\!3847 & \begin{array}{c}\langle M_1, M_2, M_3, 2\specelt\,\rangle \end{array} & 227.021533461314 & 1\cr 
3, 4, 7 & 2^{10}\!\cdot\!3^3\!\cdot\!5^3\!\cdot\!7 & \begin{array}{c}\langle M_1, M_2, M_3, 2\specelt\,\rangle \end{array} & 242.380339993035 & 1\cr 
2, 3, 8 & 2^6\!\cdot\!3^3\!\cdot\!5\!\cdot\!29\!\cdot\!103 & \begin{array}{c}\langle M_1, M_2, M_3, 2\specelt\,\rangle \end{array} & 279.394212255531 & 1\cr 
2, 3, 9 & 2^6\!\cdot\!3^2\!\cdot\!5^3\!\cdot\!7\!\cdot\!61 & \begin{array}{c}\langle M_1, M_2, M_3, 2\specelt\,\rangle \end{array} & 357.850908392946 & 2\cr 
2, 4, 8 & 2^{11}\!\cdot\!3\!\cdot\!5^3\!\cdot\!53 & \begin{array}{c}\langle M_1, M_2, M_3, 2\specelt\,\rangle \end{array} & 551.668430135999 & 2\cr 
3, 5, 6 & 2^6\!\cdot\!3^2\!\cdot\!5\!\cdot\!21737 & \begin{array}{c}\langle M_1, M_2, M_3, 2\specelt\,\rangle \end{array} & 756.149101085614 & 2\!\cdot\!3\cr 
1, 4, 7 & 2^7\!\cdot\!3^2\!\cdot\!5^3\!\cdot\!7\!\cdot\!109 & \begin{array}{c}\langle M_1, M_2, M_3, 2\specelt\,\rangle \end{array} & 1246.076904950150 & 2\!\cdot\!3\cr 
2, 6, 8 & 2^{15}\!\cdot\!3^3\!\cdot\!5^3 & \begin{array}{c}\langle M_1, M_2, M_3, 2\specelt\,\rangle \end{array} & 1351.447382022975 & 2^2\cr 
1, 4, 6 & 2^6\!\cdot\!3^3\!\cdot\!5\!\cdot\!18679 & \begin{array}{c}\langle M_1, M_2, M_3, 2\specelt\,\rangle \end{array} & 1854.823840399464 & 2^2\!\cdot\!3\cr 
1, 2, 7 & 2^6\!\cdot\!3\!\cdot\!5\!\cdot\!7\!\cdot\!47\!\cdot\!569 & \begin{array}{c}\langle M_1, M_2, M_3, 2\specelt\,\rangle \end{array} & 1970.441147103680 & 3\!\cdot\!5\cr 
1, 3, 9 & 2^6\!\cdot\!3^2\!\cdot\!5^3\!\cdot\!4273 & \begin{array}{c}\langle M_1, M_2, M_3, 2\specelt\,\rangle \end{array} & 3214.588641563687 & 2\!\cdot\!7\cr 
\exendtable
\end{table}

\begin{table}[h]
\ \bigskip\ \bigskip
\caption{Genus 4 curves $y^2+(2x^5\!+\!(v_1x\!+\!1)\cdots(v_4x\!+\!1))y+x^{10}=0\>$}
\label{tru4}
\extable{v}
1, 2, 3, 4 & 2^9\!\cdot\!3^2\!\cdot\!17\!\cdot\!113 & \begin{array}{c}\langle M_1, M_2, M_3, M_4,2\specelt\,\rangle \end{array} & 2.497987723694 & 1/2^5\!\cdot\!5\cr 
2, 3, 4, 6 & 2^8\!\cdot\!3^2\!\cdot\!80251 & \begin{array}{c}\langle M_1, M_2, M_3, M_4,2\specelt\,\rangle \end{array} & 45.006091920102 & 1/2^5\cr 
1, 2, 3, 5 & 2^9\!\cdot\!3^2\!\cdot\!5\!\cdot\!10007 & \begin{array}{c}\langle M_1, M_2, M_3, M_4,2\specelt\,\rangle \end{array} & 68.192003860287 & 1/2^3\cr 
\exendtable
\end{table}

\clearpage

\bibliographystyle{plain}
\bibliography{article}

\end{document}